\documentclass[a4paper,
reqno,
12pt,
titlepage,
tbtags
]{amsart} %

\usepackage{amsfonts}
\usepackage{amssymb}
\usepackage{mathrsfs}
\usepackage{amsmath}
\usepackage{mathtools}
\usepackage{xcolor}
\usepackage{dsfont} 
\usepackage{paralist}
\usepackage{microtype}

\usepackage[pagewise]{lineno}

\usepackage{hyperref}
\hypersetup{urlcolor = red, citecolor = blue}
\hypersetup{
  pdfauthor = Xuan Mao,
  pdftitle = {Logarithmic energy-dissipation estimates and finite-time blow-up for chemotaxis systems in critical dimensions},
  pdfsubject = Manuscript,
  pdfkeywords = {chemotaxis; logarithmic energy-dissipation estimates; critical entropy; cumulative mass; critical dimension; finite-time blow-up}
}

\hypersetup{%
	unicode=true,%
	colorlinks=true,%
	linkcolor=blue,%
	urlcolor=black,%
	citecolor=pink,%
	bookmarksopen=true,%
	bookmarksnumbered,%
	pdfpagemode=UseOutlines,%
	pdfstartview=Fit,%
	pdfpagelayout=OneColumn
  }

\usepackage[top = 1in,bottom = 1in,left = 1in,right = 1in]{geometry}
\usepackage[sort]{cite}
\usepackage[capitalize,nameinlink,noabbrev]{cleveref}

\allowdisplaybreaks[4]
\theoremstyle{plain}
\newtheorem{theorem}{Theorem}[section]

\newtheorem{lemma}[theorem]{Lemma}
\newtheorem{proposition}[theorem]{Proposition}

\theoremstyle{definition}

\theoremstyle{remark}
\newtheorem{remark}[theorem]{Remark}

\numberwithin{equation}{section}

\AddToHook{env/corollary/begin}{\crefalias{theorem}{corollary}}
\AddToHook{env/lemma/begin}{\crefalias{theorem}{lemma}}
\AddToHook{env/proposition/begin}{\crefalias{theorem}{proposition}}
\AddToHook{env/conjecture/begin}{\crefalias{theorem}{conjecture}}
\AddToHook{env/remark/begin}{\crefalias{theorem}{remark}}

\crefname{theorem}{Theorem}{Theorems}
\Crefname{theorem}{Theorem}{Theorems}
\crefname{corollary}{Corollary}{Corollaries}
\Crefname{corollary}{Corollary}{Corollaries}
\crefname{lemma}{Lemma}{Lemmas}
\Crefname{lemma}{Lemma}{Lemmas}
\crefname{proposition}{Proposition}{Propositions}
\Crefname{proposition}{Proposition}{Propositions}
\crefname{conjecture}{Conjecture}{Conjectures}
\Crefname{conjecture}{Conjecture}{Conjectures}
\crefname{definition}{Definition}{Definitions}
\Crefname{definition}{Definition}{Definitions}
\crefname{example}{Example}{Examples}
\Crefname{example}{Example}{Examples}
\crefname{remark}{Remark}{Remarks}
\Crefname{remark}{Remark}{Remarks}

\newcommand{\dd}{\;\mathrm{d}}

\providecommand{\e}{\ensuremath{\mathrm{e}}}

\begin{document}

\title[Logarithmic energy--dissipation estimates]{Logarithmic energy--dissipation estimates for radial chemotaxis systems in critical dimensions}

\author[X.~Mao]{Xuan Mao}
\address[X.~Mao]{School of Mathematics, Hohai University, Nanjing 211100, Jiangsu, China}
\address[X.~Mao]{Laboratory of Mathematical Modeling and Intelligent Computing for Water Systems, Hohai University, Nanjing 211100, China}
\email[X.~Mao]{maox@hhu.edu.cn}
\urladdr[X.~Mao]{https://orcid.org/0009-0001-4960-8901}
\thanks{Xuan~Mao is funded by
  ``the Fundamental Research Funds for the Central Universities'' (No.~B250201215), Basic Research Program of Jiangsu~(No.~BK20251482) and National Natural Science Foundation of China (No.~12601368).}


\author[M.~Liu]{Meng Liu}
\address[M.~Liu]{Department of Applied Mathematics, Anhui University of Technology, Ma'anshan
243002, P. R. China}
\email[M.~Liu]{LMeng\_math@yeah.net}
\thanks{Meng Liu is supported by the Youth Research Fund of Anhui University of Technology (No. QZ202421, No. QD202377 ).}

\author[Y.~Li]{Yuxiang Li}
\address{School of Mathematics, Southeast University, Nanjing 211189, P. R. China}
\thanks{Yuxiang Li is supported in part by National Natural Science Foundation of China (No. 12271092, No. 11671079).}
\email{lieyx@seu.edu.cn}



\begin{abstract}
  Inspired by Mizoguchi--Winkler (2014),
  we improve their sublinear estimates to a logarithmic energy--dissipation functional inequality associated with the
  parabolic--parabolic Keller--Segel system,
  valid for suitable radial function pairs in a disk.

  Our main application concerns the fully parabolic indirect
  chemotaxis system
  \begin{equation*}
    \begin{cases}
      u_t=\Delta u-\nabla\cdot(u\nabla v),\\
      v_t=\Delta v-v+w,\\
      w_t=\Delta w-w+u
    \end{cases}
  \end{equation*}
  with homogeneous Neumann boundary conditions in a ball
  $B_R\subset\mathbb R^4$.
 A closed logarithmic spacetime energy--dissipation inequality is valid along the radial solution trajectory.
  Consequently, for every prescribed supercritical mass
  $m>64\pi^2$, there exist smooth positive radially symmetric
  initial data with cell mass $m$ such that the corresponding classical solution blows
  up in finite time.
\end{abstract}

\maketitle

\section{Introduction}\label{introduce section}

We concern finite-time singularity formation for the fully
parabolic chemotaxis system with indirect signal production
\begin{equation}
  \label{sys: fpisp}
  \begin{cases}
    \partial_t u=\Delta u-\nabla\cdot(u\nabla v)
      &\text{for }x\in\Omega,\ t>0,\\
    \partial_t v=\Delta v-v+w
      &\text{for }x\in\Omega,\ t>0,\\
    \partial_t w=\Delta w-w+u
      &\text{for }x\in\Omega,\ t>0,\\
    \partial_\nu u=\partial_\nu v=\partial_\nu w=0
      &\text{for }x\in\partial\Omega,\ t>0,\\
    (u,v,w)(x,0)=(u_0,v_0,w_0)(x)
      &\text{for }x\in\Omega.
  \end{cases}
\end{equation}
in a ball $\Omega=B_R\subset\mathbb R^4$ for some $R > 0$,
where $\partial_\nu$ denotes the derivative with respect to the outward normal of $\partial\Omega$,
and initial functions $u_0\in C^0(\overline\Omega)$ and
$v_0, w_0\in C^2(\overline\Omega)$ are nonnegative and radially symmetric with $\partial_\nu v_0 = \partial_\nu w_0 = 0$ on $\partial\Omega$.

Indirect signal production separates the population that responds to
chemotactic cues from the intermediate population or substance that
produces them. Related mechanisms arise in models of beetle aggregation
\cite{Strohm2013} and phenotypic heterogeneity \cite{Macfarlane2022}.
The two successive parabolic signal equations in \eqref{sys: fpisp}
also connect the system to diffusion with distributed memory
\cite{Shi2021}.
See \cite{Winkler2025} for a broader discussion.

We first recall the classical Keller--Segel
system with direct signal production
\begin{equation}
  \label{eq:2d-ks-system}
  \begin{cases}
    u_t=\Delta u-\nabla\cdot(u\nabla v),\\
    v_t=\Delta v-v+u.
  \end{cases}
\end{equation}
Write $m=\int_\Omega u_0$ for the cell mass.
Solutions are globally bounded in one dimension \cite{Osaki2001}.
Boundedness holds for $m<4\pi$ on general smooth
bounded planar domains and for $m<8\pi$ for radial solutions in a disk
\cite{Nagai1997}. Herrero--Vel\'azquez \cite{Herrero1997}
constructed radial finite-time collapse, while Horstmann--Wang
\cite{Horstmann2001} obtained nonradial supercritical solutions
that become unbounded in finite or infinite time.
In higher dimensions, sufficiently small data in critical norms yield boundedness~\cite{Cao2015,Jiang2020},
whereas unboundedness occurs for suitable data at every prescribed mass $m>0$
on balls \cite{Winkler2010b}.

A common tool in both boundedness and blow-up theory is the
energy--dissipation identity,
$\mathcal F_{\mathrm{KS}}'(t)=-\mathcal D_{\mathrm{KS}}(t)$
along classical solutions.
The associated energy and its dissipation rate
\cite{Nagai1997,Gajewski1998,Biler1998} are
\begin{align}
  \label{eq:2d-energy-definition}
  \mathcal F_{\mathrm{KS}}(u,v)
    &=\int_\Omega u\ln u-\int_\Omega uv
      +\frac12\int_\Omega|\nabla v|^2+\frac12\int_\Omega v^2,\\
  \label{eq:2d-dissipation-definition}
  \mathcal D_{\mathrm{KS}}(u,v)
    &=\int_\Omega|\Delta v-v+u|^2
      +\int_\Omega u|\nabla(\ln u-v)|^2.
\end{align}
Nagai--Senba--Yoshida \cite{Nagai1997} combined this identity
with the Trudinger--Moser inequality to control entropy in the
two-dimensional subcritical mass. In two and higher dimensions,
comparisons of low initial energy with bounds on stationary energies
exclude globally bounded trajectories
\cite{Horstmann2001,Winkler2010b}. Such comparisons alone do not
distinguish finite-time blow-up from infinite-time grow-up.

Winkler \cite{Winkler2013} resolved this alternative in balls of
dimension $n\ge3$ by
establishing the nonlinear functional relation
\begin{equation}
  \label{eq: functional inequality}
  -\mathcal F_{\mathrm{KS}}
    \le\int_\Omega uv+\frac{|\Omega|}{\e}
    \le C\bigl(1+\mathcal D_{\mathrm{KS}}^\theta\bigr)
    \quad\text{for some }\theta\in(0,1),
\end{equation}
for suitable positive radial pairs. The sublinear power of the
dissipation and the energy identity yield finite-time blow-up
from initial data with sufficiently low energy,
which can be constructed at any mass level $m>0$.
In a disk, Mizoguchi--Winkler
\cite[Lemma~7.1]{Mizoguchi2014} obtained, combining with the radial Sobolev embedding,
\begin{equation}
  \label{eq:intro-mizoguchi-winkler-energy}
  -\mathcal F_{\mathrm{KS}}(u,v)
  \le C_q\bigl(1+m^3+\|v\|_{W^{1,q}(\Omega)}^3\bigr)
           \sqrt{1+\mathcal D_{\mathrm{KS}}(u,v)},
\end{equation}
for suitable positive radial pairs and $q\in(1,2)$.
Uniform signal estimates and the energy identity again yield a
energy-based criterion on finite-time blow-up, valid at supercritical mass level $m>8\pi$.
See \cite{Mizoguchi2020} for related finite-time blow-up criteria based on the initial energy.

For chemotaxis with indirect signal production,
the indirect effect depends on the evolution of both
signal components. In two dimensions, Tao--Winkler \cite{Tao2017}
proved global existence and a critical mass phenomenon for
infinite-time aggregation in a radial parabolic--elliptic--ODE
model.
Stationary-energy comparisons also underlie unboundedness results
for indirect chemotaxis \cite{Laurencot2019,Fujie2019,Mao2024a}.
Lauren\c{c}ot \cite{Laurencot2019} established an analogous
boundedness--unboundedness distinction for a parabolic--parabolic--ODE
model on general planar domains using a Lyapunov functional.
The Lyapunov structure of the fully parabolic system \eqref{sys: fpisp} was identified by
Fujie--Senba \cite{Fujie2017}.
Its energy and dissipation are
\begin{equation}
  \label{sym: energy and dissipation rate}
  \begin{aligned}
    \mathcal F(t)
      &=\int_\Omega u\ln u-\int_\Omega uv
        +\frac12\int_\Omega|(-\Delta+1)v|^2
        +\frac12\int_\Omega v_t^2,\\
    \mathcal D(t)
      &=\int_\Omega u|\nabla(\ln u-v)|^2
        +2\int_\Omega\bigl(|\nabla v_t|^2+v_t^2\bigr)
  \end{aligned}
\end{equation}
and along classical solutions they satisfy
\begin{equation}
  \label{eq: dfcdt}
  \mathcal F'(t)=-\mathcal D(t).
\end{equation}
 Fujie--Senba
\cite{Fujie2017} established boundedness in physical dimensions and
identified the four-dimensional threshold $64\pi^2$ through an
Adams type inequality. Their subsequent work
\cite{Fujie2019} constructed unbounded solutions with supercritical mass under mixed boundary conditions,
using the energy and quantization properties of stationary solutions.

Adapting the functional inequality approach \cite{Winkler2013,Mao2025},
Mao--Liu--Li \cite{Mao2025b} obtained finite-time blow-up
for the system \eqref{sys: fpisp} at
every prescribed positive mass on balls
$\Omega=B_R\subset\mathbb R^n$ with $n\ge5$,
via spacetime sublinear energy-dissipation estimates
\begin{equation*}
  \int_0^s\int_\Omega uv
  \le C(1+s)\biggl(1+\int_0^s\mathcal D(t)\dd t
  +\frac12\int_\Omega|(-\Delta+1)v(s)|^2
        +\frac12\int_\Omega v_t^2(s)\biggr)^\theta
\end{equation*}
for some $\theta\in(0,1)$,
along radial classical solutions.

Building upon comparison methods of Tao--Winkler
\cite{Tao2025} developing for cooperative parabolic systems, satisfied by cumulative radial densities to the parabolic--elliptic--parabolic variant of \eqref{sys: fpisp},
Mao--Liu--Li \cite{Mao2025a}
identified the four-dimensional critical mass $64\pi^2$ concerning finite-time chemotactic collapse on the
unit ball.

These results directly motivates our
study of finite-time blow-up for the critical
fully parabolic system \eqref{sys: fpisp} by seeking a spacetime energy inequality.

\subsection*{A two-dimensional energy--dissipation inequality}

The functional inequality \eqref{eq:intro-mizoguchi-winkler-energy} of the classical Keller--Segel system on a disk
$\Omega=B_R\subset\mathbb R^2$ provides a case study for spatial estimates.
Mizoguchi--Winkler \cite{Mizoguchi2014} introduced a core $B_\rho$ containing the
critical mass $8\pi$ to control the negative energy by the Moser--Trudinger inequality which requires zero boundary trace
and is therefore applied with cutoff errors involving $v(\rho)$, and the core radius $\rho$.
After canceling the signal energy terms and cutoff errors,
the negative energy on the surrounding annuli is essentially dominated by the auxiliary quantity
\begin{equation}
  \label{eq:intro-2d-excess-moment}
  Q_R[u]\coloneq\int_0^R\frac{(M(r)-8\pi)_+}{r}\dd r.
\end{equation}
Here $M(r)=\int_{B_r}u$ is the local mass, and
$(M(r)-8\pi)_+$ is its excess above the critical value $8\pi$ (see also \cref{rem:2d-local-mass-excess}).
Adapting the accumulated mass argument of
\cite[Lemma~5.5]{Mizoguchi2014}, we obtain the pointwise bound
\eqref{eq:2d-mass-excess-pointwise} for the mass excess and the
logarithmic bound \eqref{eq:2d-excess-moment-bound} for $Q_R[u]$.

To simplify and refine their critical-core strategy, we organize our argument
around the reduced energy \eqref{eq:2d-reduced-energy} in terms of cumulative
mass,
which leads us to a different functional tool, the critical
entropy--cumulative mass inequality \eqref{eq:intro-radial-entropy-choice}.
The latter applies directly on $B_R$ to the sharp truncation $h=u\mathbf1_{B_{r_*}}$, where
$M(r_*)=\min\{8\pi,m\}$, without cutoff errors. The logarithmic bound
for $Q_R[u]$ controls the remaining negative contribution.

\begin{proposition}
\label{prop:2d-sobolev-energy-dissipation}
Let $m,R>0$, and $\Omega=B_R\subset\mathbb R^2$.
Define
\begin{equation}
\label{eq:2d-admissible-pairs}
\mathcal S(m,R)\coloneq
\left\{(u,v)\;\middle|\;
\begin{aligned}
&u\in C^1(\bar\Omega),\quad v\in C^2(\bar\Omega),\quad u,v\text{ are radially symmetric},\\
&u,v>0\quad\text{on }\bar\Omega,\quad
\textstyle\int_\Omega u=m,\quad \partial_\nu v=0\quad\text{on }\partial\Omega
\end{aligned}
\right\}.
\end{equation}
There exists $C_R>0$, depending only on $R$, such that every
$(u,v)\in\mathcal S(m,R)$ satisfies
\begin{equation}
\label{eq:2d-sobolev-energy-dissipation}
-\frac{\mathcal F_{\mathrm{KS}}(u,v)}{C_R(1+m^2)}\le
\begin{cases}
  1 & \text{if }0<m\le8\pi,\\
  1+\ln(1+\mathcal D_{\mathrm{KS}}(u,v)) & \text{if }m>8\pi.
\end{cases}
\end{equation}
\end{proposition}

\begin{remark}
The sharpness for the logarithmic growth is provided in \cref{sec:2d-supercritical-data}.
In higher dimensions, $\theta=\frac{n-2}{n}$ can be expected to be optimal in \eqref{eq: functional inequality}.
\end{remark}

\subsection*{Main results}

We now turn to the fully parabolic system \eqref{sys: fpisp} on
$\Omega=B_R\subset\mathbb R^4$.
Building on the critical-core strategy of Mizoguchi--Winkler
\cite{Mizoguchi2014}, the energy--dissipation inequality of
\cite{Winkler2013}, and the spacetime viewpoint of
\cite{Mao2025b},
as well as energy reduction discussed above,
we establish the following logarithmic energy--dissipation estimate.

\begin{proposition}
  \label{prop:4d-closed-spacetime}
  Let $m,R,L>0$ and $\Omega=B_R\subset\mathbb R^4$.
  Define
\begin{equation}
  \label{h: initial data}
  \mathcal I(m,R)\coloneq
  \left\{
  \begin{aligned}
    &(u_0,v_0,w_0)\colon
  u_0\in C(\bar\Omega), v_0,w_0\in C^2(\bar\Omega) \text{ are radially symmetric}\\
  &\text{and nonnegative}
  \text{ such that }\textstyle\int_\Omega u_0=m,
  \partial_\nu v_0=\partial_\nu w_0=0\text{ on }\partial\Omega.
  \end{aligned}
  \right\}.
\end{equation}
  There exists $C_R>0$, depending only on $R$, such that
  the classical solution of \eqref{sys: fpisp} obeys
\begin{equation}
  \label{eq:4d-intro-closed}
  \begin{aligned}
    \int_0^s(-\mathcal F(t))_+\dd t
    &\le C_R(1+m^2)\ln\left(e+\frac{L^2}{1+m^2}\right)s\\
    &\quad\times\left[1+\ln\left(1+
       \frac{1+\int_0^s\mathcal D(t)\dd t}{s}\right)\right]
       \quad\text{for }s\in(0,T_{\max}),
  \end{aligned}
\end{equation}
whenever
  $(u_0,v_0,w_0)\in\mathcal I(m,R)$ satisfies
  \begin{equation}
    \label{eq:4d-initial-weak-bound}
    \|w_0\|_{W^{1,5/4}(\Omega)}
         +\|v_0\|_{W^{2,7/4}(\Omega)}\le L.
  \end{equation}
Here $a_+=\max\{a,0\}$, and $\mathcal F$ and $\mathcal D$ are
given in \eqref{sym: energy and dissipation rate}.
\end{proposition}

Combined with the energy identity, the logarithmic estimate gives
an initial energy criterion for finite-time blow-up.
The concentrating data of Fujie--Senba \cite[Section~4]{Fujie2019}
allow us to apply the criterion at every supercritical mass.

\begin{theorem}
  \label{thm:4d-low-energy-blowup}
  Let $R>0$ and $\Omega=B_R\subset\mathbb R^4$.
  Write $\mathcal F(0)$ the initial energy \eqref{eq:4d-initial-energy}.
  \begin{enumerate}[(i)]
    \item For every $m,L>0$, there exist $K=K(m,R,L)>0$ and
    $T=T(m,R,L)>0$ such that, for any initial data
    $(u_0,v_0,w_0)\in\mathcal I(m,R)$ satisfying
    \eqref{eq:4d-initial-weak-bound} and
    \begin{equation}
      \label{eq:4d-low-energy-data}
      \mathcal F(0)\le-K,
    \end{equation}
    the corresponding solution of \eqref{sys: fpisp} blows up
    in finite time, with $T_{\max}\le T$.

    \item
    For every $m>64\pi^2$, there exist smooth, strictly positive
    radial initial data in $\mathcal I(m,R)$ for which the classical
    solution of \eqref{sys: fpisp} blows up in finite time.
  \end{enumerate}
\end{theorem}

\begin{remark}
  The uniform integrability theory in \cite{Bai2018} gives
  $\limsup_{t\to T_{\max}}\|u\|_{p}=\infty$ for $p\in(1,\infty]$ if $T_{\max} < \infty$.
  The radial solutions of the system \eqref{sys: fpisp} exist globally and remain bounded~\cite{Fujie2017} if $m<64\pi^2$.
\end{remark}

\cref{sec:radial-entropy-proof} states and proves the critical
entropy--cumulative mass inequality and explains its relation
to the classical logarithmic Hardy--Littlewood--Sobolev inequality.
\cref{sec:two-dimensional-supplement} proves
\cref{prop:2d-sobolev-energy-dissipation} using this inequality
and local mass estimates.
\cref{sec:4d-proof} proves the four-dimensional energy--dissipation
estimate and the finite-time blow-up results. The spatial energy
and local mass estimates in \cref{sec:4d-spatial-proof} follow the
structure of the two-dimensional case study.
\cref{sec:4d-closed-proof} establishes time-integrated control
of the elliptic residual, and \cref{sec:4d-blowup-proof} combines
the spatial and temporal estimates to give the closed spacetime
bound, then proves the low-energy blow-up criterion and applies
it to the concentrating data of Fujie--Senba.
\cref{sec:2d-supercritical-data} shows the optimality of the
two-dimensional logarithmic bound by constructing concentrating
data.

\section{The critical entropy--cumulative mass inequality}
\label{sec:radial-entropy-proof}

The entropy estimate used in both critical dimensions follows from a
one-dimensional inequality for probability densities on $(0,1)$
through the radial volume variable.

\begin{lemma}
\label{lem:one-dimensional-entropy}
Let $q\in L^\infty((0,1))$ be nonnegative with $\int_0^1q(s)\dd s=1$,
and set
\begin{equation*}
Q(s)=\int_0^s q(\sigma)\dd\sigma
\quad\text{for }s\in[0,1].
\end{equation*}
With $0\ln0=0$, one has
\begin{equation}
\label{eq:one-dimensional-entropy}
J[q]\coloneq\int_0^1q(s)\ln q(s)\dd s
-\int_0^1\frac{Q(s)^2}{s}\dd s\ge-1.
\end{equation}
\end{lemma}

\begin{proof}
The function $Q$ is absolutely continuous, with $Q'=q$ almost
everywhere, $Q(0)=0$, and $Q(1)=1$. Since
$0\le Q(s)\le\|q\|_\infty s$, the endpoint $Q(s)\ln s$ tends to zero
as $s\searrow0$. Integration by parts therefore gives
\begin{equation*}
\int_0^1q(s)\ln s\dd s=-\int_0^1\frac{Q(s)}s\dd s,
\end{equation*}
and hence
\begin{equation*}
J[q]=\int_0^1\left[q(s)\ln\big(sq(s)\big)
  +\frac{Q(s)(1-Q(s))}{s}\right]\dd s.
\end{equation*}
The relative-entropy inequality reads
\begin{equation*}
a\ln a+b\ge a\ln b+a
\quad\text{for }a\ge0\text{ and }b>0.
\end{equation*}
Fix $\delta\in(0,1/2)$ and let
$E_\delta=\{s\in(0,1):\delta<Q(s)<1-\delta\}$.
Since $Q$ is continuous and nondecreasing, $E_\delta$ is an interval
with endpoint values $\delta$ and $1-\delta$.
Apply the relative-entropy inequality with $a=sq(s)$ and
$b=Q(s)(1-Q(s))$, and divide by $s$. The absolute-continuity
chain rule for a primitive of $z\mapsto\ln(z(1-z))$ gives
\begin{align*}
&\int_{E_\delta}\left[q(s)\ln\big(sq(s)\big)
  +\frac{Q(s)(1-Q(s))}{s}\right]\dd s\\
&\quad\ge\int_{E_\delta}q(s)\ln\big(Q(s)(1-Q(s))\big)\dd s
  +\int_{E_\delta}q(s)\dd s\\
&\quad=\int_\delta^{1-\delta}\ln\big(z(1-z)\big)\dd z+1-2\delta
\quad\text{for }\delta\in(0,1/2).
\end{align*}
The integrand on the left is integrable because $q\ln q$,
$q\ln s$, and $Q/s$ are integrable. On $\{Q=0\}\cup\{Q=1\}$,
one has $q=0$ almost everywhere. Dominated convergence on the left
and integrability of the logarithmic endpoints on the right yield
\begin{equation*}
J[q]\ge\int_0^1\ln\big(z(1-z)\big)\dd z+1=-1
\end{equation*}
as $\delta\searrow0$.
\end{proof}

\begin{remark}
The lower bound in \cref{prop:radial-entropy} is sharp and is not
attained. Indeed, for $a>0$ define the probability density
\begin{equation*}
q_a(s)=\frac{a(1+a)}{(a+s)^2},\qquad
Q_a(s)=\frac{(1+a)s}{a+s}
\quad\text{for }s\in[0,1].
\end{equation*}
A direct calculation gives
\begin{equation*}
J[q_a]=-1+a-a^2\ln(1+1/a)\to-1
\quad\text{as }a\searrow0.
\end{equation*}
On the other hand, equality in the relative-entropy step of the
proof of \cref{lem:one-dimensional-entropy} would require
$sq(s)=Q(s)(1-Q(s))$ almost everywhere on $\{0<Q<1\}$.
Since $q=0$ almost everywhere on $\{Q=0\}\cup\{Q=1\}$,
we would have
\begin{equation*}
sQ'(s)=Q(s)(1-Q(s))
\quad\text{for almost every }s\in(0,1).
\end{equation*}
For every $\varepsilon\in(0,1)$, the terminal condition $Q(1)=1$
and uniqueness for this scalar equation on $[\varepsilon,1]$ would
give $Q\equiv1$ there, contradicting $Q(0)=0$ and continuity.
\end{remark}

\begin{proposition}[Critical entropy--cumulative mass inequality]
\label{prop:radial-entropy}
Let $n\ge2$ be an integer and let $R,\mu>0$. Suppose that
$h\in L^\infty(B_R)$ is nonnegative and radial on
$B_R\subset\mathbb R^n$, with $\int_{B_R}h=\mu$, and set
$M_h(r)=\int_{B_r}h$. With $0\ln0=0$, one has
\begin{equation}
  \label{eq:intro-radial-entropy-choice}
  \int_{B_R}h\ln h
      -\frac n\mu\int_0^R\frac{M_h(r)^2}{r}\dd r
  \ge\mu\left(\ln\frac{\mu}{|B_R|}-1\right).
\end{equation}
\end{proposition}

\begin{proof}
Write $\omega_n=|\mathbb S^{n-1}|$, so that
$|B_R|=\omega_nR^n/n$. Introduce the radial volume variable and
the normalized density
\begin{equation*}
s=(r/R)^n,\qquad
q(s)=\frac{|B_R|}{\mu}h(Rs^{1/n})
\quad\text{for }s\in(0,1).
\end{equation*}
Then $q\in L^\infty((0,1))$ is nonnegative with $\int_0^1q=1$.
Its cumulative function $Q$ satisfies
\begin{equation*}
M_h(r)=\mu Q\big((r/R)^n\big)
\quad\text{for }r\in[0,R].
\end{equation*}
For radial integrands, the volume element becomes $|B_R|\dd s$.
Using $\dd r/r=\dd s/(ns)$, we obtain
\begin{equation}
\label{eq:entropy-cumulative-mass-volume}
\frac1\mu\int_{B_R}h\ln h
-\frac n{\mu^2}\int_0^R\frac{M_h(r)^2}{r}\dd r
=\ln\frac\mu{|B_R|}+J[q].
\end{equation}
Combining \eqref{eq:entropy-cumulative-mass-volume} with
\eqref{eq:one-dimensional-entropy} gives \eqref{eq:intro-radial-entropy-choice}.
\end{proof}

\begin{remark}
\label{rem:radial-entropy-log-hls}
In dimension two,
\eqref{eq:intro-radial-entropy-choice} follows from
the classical sharp logarithmic Hardy--Littlewood--Sobolev inequality
\cite{Carlen1992,Carlen2010},
which reads
\begin{equation*}
\int_{\mathbb R^2}h\ln h
+\frac2\mu\iint_{\mathbb R^2\times\mathbb R^2}
 h(x)h(y)\ln|x-y|\dd x\dd y
\ge\mu\left(\ln\frac\mu\pi-1\right),
\end{equation*}
valid for $h$ extended by zero outside $B_R$.
The angular
identity \cite[Section~2.3]{Calvez2012}
\begin{equation*}
\frac1{2\pi}\int_0^{2\pi}\ln|re^{i\theta}-s|\dd\theta
=\ln\max\{r,s\}
\quad\text{for }r,s>0
\end{equation*}
and integration by parts imply
\begin{equation*}
\iint_{\mathbb R^2\times\mathbb R^2}
 h(x)h(y)\ln|x-y|\dd x\dd y
=\mu^2\ln R-\int_0^R\frac{M_h(r)^2}{r}\dd r,
\end{equation*}
which indeed reduces HLS inequality into \eqref{eq:intro-radial-entropy-choice} for $n=2$ with the stated constant.
\end{remark}

\section{A case study: the two-dimensional Keller--Segel system}
\label{sec:two-dimensional-supplement}

This section is devoted to the logarithmic energy--dissipation estimate for the
classical parabolic--parabolic Keller--Segel system
\eqref{eq:2d-ks-system} in a disk $\Omega=B_R\subset\mathbb R^2$,
via the critical entropy--cumulative mass inequality and local mass
estimates.

Fix $m>0$ and $R>0$. Unless otherwise specified, let
$(u,v)\in\mathcal S(m,R)$, with the class defined in
\eqref{eq:2d-admissible-pairs}.
Following Winkler \cite{Winkler2013} and Mizoguchi--Winkler~\cite{Mizoguchi2014}, we introduce $f$ and $g$ to denote the quantities in dissipation accordingly
\begin{equation}
  \label{eq:2d-dissipation-variables}
  \begin{cases}
    -\Delta v+v=u+f & \text{in }\Omega,\\
    u_r - uv_r = \sqrt u\,g & \text{for }r\in(0,R],\\
    \partial_\nu v=0 & \text{on }\partial\Omega.
  \end{cases}
\end{equation}
Then \eqref{eq:2d-dissipation-definition} gives
\begin{equation}
  \label{eq:2d-dissipation-decomposition}
  \mathcal D_{\mathrm{KS}}(u,v)=\|f\|_2^2+\|g\|_2^2.
\end{equation}
Define the local mass function and the accumulated signal error by
\begin{equation}
\label{eq:2d-mass-definitions}
M(r)=2\pi\int_0^r s u(s)\dd s,\qquad
E(r)=2\pi\int_0^r s\big(v(s)-f(s)\big)\dd s
\quad\text{for }r\in[0,R].
\end{equation}
Thus integration of
$(rv_r)'=r(v-f-u)$ from zero to $r$ gives the pointwise representation
\begin{equation}
\label{eq:2d-radial-signal}
rv_r(r)=\frac{E(r)-M(r)}{2\pi}
\quad\text{for }r\in(0,R].
\end{equation}
Here and below, $C_R>0$ denotes a constant depending
only on $R$, whose value may change from line to line.

We begin with the radial energy representation.
Combining the interaction term with the gradient part of the signal
energy cancels the mixed signal contributions exactly.
Related cancellations are more transparent in parabolic--elliptic
chemotaxis systems \cite{Carrillo2014, Mao2026}: the elliptic signal equation combines the
interaction term with the quadratic signal energy, reducing the
free energy to entropy and a nonlocal interaction of the density.

\begin{lemma}
\label{lem:2d-energy-representation}
There exists $C_R>0$ such that
\begin{equation}
\label{eq:2d-reduced-energy}
\mathcal F_{\mathrm{KS}}(u,v)
\ge\int_\Omega u\ln u
-\frac1{4\pi}\int_0^R\frac{M(r)^2}{r}\dd r-C_Rm^2
\quad\text{for } (u,v)\in\mathcal S(m,R).
\end{equation}
\end{lemma}

\begin{proof}
Integration by parts in $\int_0^RvM'$ and
\eqref{eq:2d-radial-signal} yield
\begin{align*}
\int_\Omega uv
&=mv(R)-\int_0^RMv_r\dd r
=mv(R)+\frac1{2\pi}\int_0^R\frac{M(M-E)}r\dd r,\\
\frac12\int_\Omega|\nabla v|^2
&=\frac1{4\pi}\int_0^R\frac{(M-E)^2}r\dd r.
\end{align*}
Substituting these identities into
\eqref{eq:2d-energy-definition} cancels the mixed $ME$ terms
and gives
\begin{equation}
\label{eq:2d-exact-energy-representation}
\mathcal F_{\mathrm{KS}}(u,v)
=\int_\Omega u\ln u
-\frac1{4\pi}\int_0^R\frac{M^2}r\dd r
+\frac1{4\pi}\int_0^R\frac{E^2}r\dd r
+\frac12\int_\Omega v^2-mv(R).
\end{equation}

To bound the boundary term, integrate $(r^2v(r))'$ and use
\eqref{eq:2d-radial-signal} to obtain
\begin{equation*}
v(R)=\frac1{\pi R^2}\int_\Omega v
 +\frac1{2\pi R^2}\int_0^R r\bigl(E(r)-M(r)\bigr)\dd r.
\end{equation*}
Since $M\ge0$, Cauchy--Schwarz and
$\int_0^Rr^3\dd r=R^4/4$ imply
\begin{equation*}
v(R)\le\frac{\|v\|_2}{\sqrt\pi R}
 +\frac1{4\pi}\left(\int_0^R\frac{E(r)^2}{r}\dd r\right)^{1/2}.
\end{equation*}
Young's inequality therefore yields
\begin{equation}
\label{eq:2d-boundary-absorption}
mv(R)\le\frac12\|v\|_2^2
 +\frac1{4\pi}\int_0^R\frac{E(r)^2}{r}\dd r+C_Rm^2,
\end{equation}
which reduces
\eqref{eq:2d-exact-energy-representation} into
 \eqref{eq:2d-reduced-energy}.
\end{proof}

The critical entropy inequality  reduces the negative
energy to the weighted integral $Q_R[u]$ of the local mass excess
above $8\pi$, defined in
\eqref{eq:intro-2d-excess-moment}.

\begin{lemma}
\label{lem:2d-energy-reduction}
It holds that
\begin{equation}
\label{eq:2d-energy-reduction}
-\mathcal F_{\mathrm{KS}}(u,v)
\le C_R(1+m^2)+\frac{m+8\pi}{4\pi}Q_R[u]
\quad\text{for } (u,v)\in\mathcal S(m,R).
\end{equation}
\end{lemma}

\begin{proof}
Set $\mu=\min\{m,8\pi\}$. Since $M$ is strictly increasing from
zero to $m$, there is a unique $r_\mu\in(0,R]$ satisfying
$M(r_\mu)=\mu$.
Define $h=u\mathbf1_{B_{r_\mu}}$ and $u_{\mathrm{o}}=u\mathbf1_{\Omega\setminus B_{r_\mu}}$. Then $h$ is bounded,
nonnegative, and radial, with $\int_\Omega h=\mu$, and its local mass is
\begin{equation*}
M_h(r)=\int_{B_r}h=\min\{M(r),8\pi\}
\quad\text{for }r\in[0,R].
\end{equation*}

Applying \cref{prop:radial-entropy} with $n=2$ and using
$2/\mu\ge1/(4\pi)$ gives
\begin{equation}
\label{eq:2d-critical-entropy}
\int_\Omega h\ln h-\frac1{4\pi}\int_0^R\frac{M_h(r)^2}{r}\dd r
\ge\mu\left(\ln\frac{\mu}{\pi R^2}-1\right)\ge-\pi R^2.
\end{equation}
Moreover, $M\leq M(R)=m$ implies
\begin{align}
M(r)^2-M_h(r)^2
&=(M(r)-8\pi)_+(M(r)+8\pi)\notag\\
&\le(m+8\pi)(M(r)-8\pi)_+
\quad\text{for }r\in[0,R].
\label{eq:2d-squared-mass-excess}
\end{align}
The cutoff function $\mathbf1_{B_{r_\mu}}$ gives
an exact decomposition of the entropy over $B_{r_\mu}$ and its
complement. Adding and subtracting the squared local mass of $h$,
we obtain
\begin{align}
&\int_\Omega u\ln u-\frac1{4\pi}\int_0^R\frac{M(r)^2}{r}\dd r\notag\\
&\quad=\left[\int_\Omega h\ln h
   -\frac1{4\pi}\int_0^R\frac{M_h(r)^2}{r}\dd r\right]\notag\\
&\qquad+\int_{\Omega}u_{\mathrm{o}}\ln u_{\mathrm{o}}
   -\frac1{4\pi}\int_0^R\frac{M(r)^2-M_h(r)^2}{r}\dd r\notag\\
&\quad\ge-\pi R^2-\frac{|\Omega\setminus B_{r_\mu}|}{e}
   -\frac{m+8\pi}{4\pi}Q_R[u]
\ge-2\pi R^2-\frac{m+8\pi}{4\pi}Q_R[u].
\label{eq:2d-core-excess-energy}
\end{align}
The first inequality in \eqref{eq:2d-core-excess-energy} follows
from \eqref{eq:2d-critical-entropy}, the bound
$u_{\mathrm{o}}\ln u_{\mathrm{o}}\ge-1/e$ on $\Omega\setminus B_{r_\mu}$,
and integration of \eqref{eq:2d-squared-mass-excess} against
$\dd r/r$ using \eqref{eq:intro-2d-excess-moment}.
Thus \eqref{eq:2d-reduced-energy} and \eqref{eq:2d-core-excess-energy}
imply \eqref{eq:2d-energy-reduction}.
\end{proof}

\begin{remark}
\label{rem:2d-local-mass-excess}
The cumulative mass of the outer component $u_{\mathrm{o}}$ is exactly the local
mass excess:
\begin{equation*}
\int_{B_r}u_{\mathrm{o}}
=M(r)-\min\{M(r),8\pi\}=(M(r)-8\pi)_+
\quad\text{for }r\in[0,R].
\end{equation*}
Fubini's theorem therefore gives
\begin{equation*}
Q_R[u]=\int_\Omega u_{\mathrm{o}}(x)\ln\frac{R}{|x|}\dd x.
\end{equation*}
Thus $Q_R[u]$ is the logarithmic moment of the outer component.
\end{remark}

It remains to control $Q_R[u]$ by the dissipation.
Motivated by \cite[Lemma~5.5]{Mizoguchi2014},
we estimate
the local mass excess pointwise and then integrate the resulting bound.

\begin{lemma}
\label{lem:2d-local-mass}
There exists $C_R>0$ such that whenever $(u,v)\in\mathcal S(m,R)$,
\begin{equation}
\label{eq:2d-mass-excess-pointwise}
\big(M(r)-8\pi\big)_+
\le C_Rr\big(m+\sqrt{\mathcal D_{\mathrm{KS}}(u,v)}\big)
\quad\text{for }r\in(0,R].
\end{equation}
Moreover,
\begin{equation}
\label{eq:2d-excess-moment-bound}
Q_R[u]\le C_Rm\big(1+\ln(1+\mathcal D_{\mathrm{KS}}(u,v))\big)
\quad\text{for } (u,v)\in\mathcal S(m,R).
\end{equation}
\end{lemma}

\begin{proof}
Since $u$ is bounded, $s^2u(s)\to0$ as $s\searrow0$.
Integrating $(s^2u(s))'=2su(s)+s^2u_r(s)$ and using the radial
equation in \eqref{eq:2d-dissipation-variables} and $M'(s)=2\pi su(s)$ yields
\begin{align*}
rM'(r)
&=2\pi r^2u(r)=2M(r)+2\pi\int_0^r s^2u_r(s)\dd s\\
&=2M(r)+\int_0^r sv_r(s)M'(s)\dd s
  +2\pi\int_0^r s^2\sqrt{u(s)}g(s)\dd s
\quad\text{for }r\in(0,R].
\end{align*}
By \eqref{eq:2d-radial-signal} and $M(0)=0$, we have
\begin{align*}
\int_0^r sv_r(s)M'(s)\dd s
&=\frac1{2\pi}\int_0^r E(s)M'(s)\dd s
  -\frac1{2\pi}\int_0^r M(s)M'(s)\dd s\\
&=\frac1{2\pi}\int_0^r E(s)M'(s)\dd s
  -\frac{M(r)^2}{4\pi}
\quad\text{for }r\in(0,R].
\end{align*}
Combining these identities, we obtain
\begin{equation}
\label{eq:2d-quadratic-mass-identity}
rM'(r)=2M(r)-\frac{M(r)^2}{4\pi}+\mathcal R(r)
\quad\text{for }r\in(0,R],
\end{equation}
where
\begin{equation*}
\mathcal R(r)\coloneq\frac1{2\pi}\int_0^r E(s)M'(s)\dd s
 +2\pi\int_0^r s^2\sqrt{u(s)}g(s)\dd s
\quad\text{for }r\in(0,R].
\end{equation*}

Applying the Neumann $L^1$-regularity estimate of
\cite[Lemma~23]{Brezis1973} to the elliptic equation in
\eqref{eq:2d-dissipation-variables} and using the two-dimensional
Sobolev embedding $W^{1,4/3}(\Omega)\hookrightarrow L^2(\Omega)$ gives
\begin{equation}
\label{eq:2d-signal-L2}
\|v\|_2\le C_R\|v\|_{W^{1,4/3}(\Omega)}
\le C_R\|u+f\|_1\le C_R\big(m+\|f\|_2\big).
\end{equation}
By Cauchy--Schwarz,
\begin{equation*}
|E(s)|\le\sqrt\pi\,s\|v-f\|_2
\quad\text{for }s\in[0,R],
\end{equation*}
and hence positivity of $M'(r)=2\pi r u(r)$ gives
\begin{align*}
\left|\frac1{2\pi}\int_0^r E(s)M'(s)\dd s\right|
&\le\frac{rM(r)}{2\sqrt\pi}\|v-f\|_2
\quad\text{for }r\in(0,R],\\
\left|2\pi\int_0^r s^2\sqrt{u(s)}g(s)\dd s\right|
&=\left|\int_{B_r}|x|\sqrt u\,g\right|
\le r\sqrt{M(r)}\,\|g\|_2
\quad\text{for }r\in(0,R].
\end{align*}
Together with \eqref{eq:2d-signal-L2} and the triangle inequality,
these estimates imply
\begin{equation}
\label{eq:2d-mass-remainder-dissipation}
|\mathcal R(r)|
\le C_Rr\left[M(r)\big(m+\|f\|_2\big)
   +\sqrt{M(r)}\,\|g\|_2\right]
\quad\text{for }r\in(0,R].
\end{equation}
For $r\in(0,R]$ with $M(r)>8\pi$, the fact that $M'\ge0$ and
\eqref{eq:2d-quadratic-mass-identity} give
\begin{align*}
M(r)-8\pi
&\le\frac{4\pi}{M(r)}|\mathcal R(r)|
\le C_Rr\left(m+\|f\|_2+\frac{\|g\|_2}{\sqrt{M(r)}}\right)\\
&\le C_Rr\big(m+\sqrt{\mathcal D_{\mathrm{KS}}}\big)
    \quad\text{for }r\in(0,R]\text{ with }M(r)>8\pi.
\end{align*}
The last step uses $M(r)>8\pi$ and \eqref{eq:2d-dissipation-decomposition}.
When $M(r)\le8\pi$, the positive part vanishes, so this proves
\eqref{eq:2d-mass-excess-pointwise} for all $r\in(0,R]$.

For the integral estimate, put
$\Lambda=C_R(1+\sqrt{\mathcal D_{\mathrm{KS}}})\ge1$.
A positive local mass excess requires $m>8\pi>1$, so
\eqref{eq:2d-mass-excess-pointwise} also gives
\begin{equation*}
(M(r)-8\pi)_+\le m\Lambda r
\quad\text{for }r\in(0,R].
\end{equation*}
Set $r_0=\min\{R,\Lambda^{-1}\}$. Splitting the integral at $r_0$
and using $(M-8\pi)_+\le m$, we obtain
\begin{align*}
Q_R[u]
&\le m\Lambda r_0+m\ln\frac R{r_0}\\
&\le m\big[1+\max\{0,\ln(R\Lambda)\}\big]\\
&\le C_Rm\big(1+\ln(1+\mathcal D_{\mathrm{KS}})\big).
\end{align*}
The last step uses
$1+\sqrt{\mathcal D_{\mathrm{KS}}}\le\sqrt{2(1+\mathcal D_{\mathrm{KS}})}$.
This proves \eqref{eq:2d-excess-moment-bound}.
\end{proof}

\begin{proof}[Proof of \cref{prop:2d-sobolev-energy-dissipation}]
If $m\le8\pi$, then $M(r)\le m\le8\pi$ for $r\in[0,R]$,
so $Q_R[u]=0$, and \eqref{eq:2d-energy-reduction} gives the first
case of \eqref{eq:2d-sobolev-energy-dissipation}.
For $m>8\pi$, inserting \eqref{eq:2d-excess-moment-bound} into
\eqref{eq:2d-energy-reduction} and using $m(m+8\pi)\le(1+4\pi)(1+m^2)$
gives the second case of \eqref{eq:2d-sobolev-energy-dissipation}.
\end{proof}

\section{Energy--dissipation estimates and blow-up in dimension four}
\label{sec:4d-proof}

We prove the four-dimensional spacetime estimate in
\cref{prop:4d-closed-spacetime} and the finite-time blow-up results
in \cref{thm:4d-low-energy-blowup} by combining spatial energy and
local mass estimates with time-integrated control of the elliptic
residual. We first introduce notations used throughout this section.

Fix $m>0$ and initial data $(u_0,v_0,w_0)\in\mathcal I(m,R)$,
with the class defined in \eqref{h: initial data}. The local theory in
\cite[Section~4]{Fujie2017} provides a unique
nonnegative radial classical solution of \eqref{sys: fpisp} on
a maximal interval $[0,T_{\max})$. The solution is smooth and positive for $t\in(0,T_{\max})$.

Let $A=-\Delta+1$ be the positive self-adjoint Neumann operator
on $L^2(\Omega)$, with
$D(A)=\{\psi\in W^{2,2}(\Omega)\colon\partial_\nu\psi=0\}$.
Following \cite[Section~2.6, formula~(6.9) and Definition~6.7]{Pazy1983},
we define the negative powers by
\begin{equation}
  \label{eq:4d-negative-powers}
  A^{-\alpha}\psi
  =\frac1{\Gamma(\alpha)}\int_0^\infty
       \tau^{\alpha-1}e^{-\tau A}\psi\dd\tau
  \quad\text{for }\alpha>0,\ \psi\in L^2(\Omega).
\end{equation}
Set $A^0=I$ and define $A^\alpha=(A^{-\alpha})^{-1}$ on
$D(A^\alpha)=\operatorname{Ran}(A^{-\alpha})$ for $\alpha>0$.

Since $v_0,w_0\in D(A)$, the two signal variation-of-constants
formulas give $Av(t)\to Av_0$ in $L^2(\Omega)$ as $t\searrow0$.
Hence $v_t(t)=w(t)-Av(t)$ tends to $w_0-Av_0$ in $L^2(\Omega)$.
Together with the uniform convergence of $u(t)$ to $u_0$, these limits show that
the energy in \eqref{sym: energy and dissipation rate} extends
continuously to $t=0$, with
\begin{equation}
  \label{eq:4d-initial-energy}
  \mathcal F(0)=\int_\Omega u_0\ln u_0-\int_\Omega u_0v_0
    +\frac12\|Av_0\|_2^2+\frac12\|w_0-Av_0\|_2^2.
\end{equation}
 Integrating \eqref{eq: dfcdt} from a positive time
and passing to zero gives
\begin{equation}
  \label{eq:4d-integrated-energy}
  \int_0^s\mathcal D(t)\dd t
     =\mathcal F(0)-\mathcal F(s)<\infty
     \quad\text{for }s\in(0,T_{\max}).
\end{equation}

We introduce the stationary residuals
\begin{equation}
  \label{eq:4d-residual-definitions}
  g=\frac{u_r}{\sqrt u}-\sqrt u\,v_r,\qquad
  \mathcal R=A^2v-u,
\end{equation}
which measure the departure from the zero-flux relation $u_r-uv_r=0$ and the stationary
signal relation $A^2v=u$, respectively.
Eliminating $w$ from the signal equations gives
\begin{equation}
  \label{eq:4d-damped-signal}
  v_{tt}+2Av_t+A^2v=u
  \quad\text{in }\Omega\times(0,T_{\max}).
\end{equation}
Using the residuals in \eqref{eq:4d-residual-definitions}, we rewrite
the density and signal equations as
\begin{equation}
  \label{eq:4d-residual-evolution}
  \begin{cases}
    0=u_t-r^{-3}(r^3\sqrt u\,g)_r,\\
    0=v_{tt}+2Av_t+\mathcal R.
  \end{cases}
\end{equation}
The Neumann form of $A$ gives
$\|A^{1/2}\psi\|_2^2=\|\nabla\psi\|_2^2+\|\psi\|_2^2$
for $\psi\in W^{1,2}(\Omega)$. Hence the dissipation in
\eqref{sym: energy and dissipation rate} takes the form
\begin{equation}
  \label{eq:4d-dissipation-operator}
  \mathcal D=2\|A^{1/2}v_t\|_2^2+\|g\|_2^2
  \quad\text{for }t\in(0,T_{\max}),
\end{equation}
which directly controls $g$ and $A^{1/2}v_t$.

We measure $\mathcal R$ in the negative-order norm
$\|A^{-3/4}\mathcal R\|_2$.
The exponent $3/4$ is chosen to
accommodate both the spatial and temporal estimates that produce the closed energy-dissipation inequality. Spatially,
this norm controls radial gradient \eqref{eq:4d-residual-gradient} of the signal error $h=v-A^{-2}u=A^{-2}\mathcal R$ pointwise.
Temporally,
integrating the differential inequality for the cross functional
$\mathcal J(t)=\langle v_t,A^{-3/2}\mathcal R\rangle$ bounds the time
integral of the squared residual norm by the accumulated dissipation \eqref{eq:4d-integrated-energy}
and the endpoint values of $\mathcal J$, as shown in
\cref{lem:4d-residual-time}.

\subsection{Spatial energy and local mass estimates}
\label{sec:4d-spatial-proof}

The logarithmic spatial estimate below
bounds the negative energy in terms of the density flux and the
elliptic residual. Throughout this subsection, we fix $m,R>0$,
 and suppress the time argument $t\in(0,T_{\max})$ in spatial quantities.

\begin{proposition}
  \label{prop:4d-static-log-energy}
  There exists $C_R>0$ such that
  the solution $(u,v,w)$ of \eqref{sys: fpisp} satisfies
  \begin{equation}
    \label{eq:4d-static-log-energy}
    (-\mathcal F)_+
      \le C_R(1+m^2)\left[1+\ln\left(
        1+\frac{\|g\|_2^2+\|A^{-3/4}\mathcal R\|_2^2}{1+m^2}
      \right)\right]
      \quad\text{for } t\in(0,T_{\max}),
  \end{equation}
  whenever $(u_0,v_0,w_0)\in\mathcal I(m,R)$.
\end{proposition}

As in \cref{sec:two-dimensional-supplement}, we reduce the negative
energy to the weighted integral of the local mass excess
\begin{equation}
  \label{eq:4d-excess-moment-definition}
  Q_R[u]\coloneq\int_0^R\frac{(M(r)-64\pi^2)_+}{r}\dd r,
  \qquad
  M(r)=\int_{B_r}u \quad\text{for }r\in[0,R].
\end{equation}

\begin{lemma}
  \label{lem:4d-energy-reduction}
  There exists $C_R>0$ such that for any
  $(u_0,v_0,w_0)\in\mathcal I(m,R)$,
  \begin{equation}
    \label{eq:4d-energy-excess}
    -\mathcal F(t)\le C_R(1+m^2)+\frac{m+64\pi^2}{16\pi^2}Q_R[u(t)]
    \quad\text{for }t\in(0,T_{\max}).
  \end{equation}
\end{lemma}

\begin{samepage}
The proof of \cref{lem:4d-energy-reduction} starts by completing the
square in the signal terms.

\begin{lemma}
  \label{lem:4d-energy-representation}
  With $\xi=A^{-1}u$, the energy satisfies
  \begin{equation}
    \label{eq:4d-energy-square}
    \mathcal F(t)
    =\int_\Omega u\ln u-\frac12\|\xi\|_2^2
       +\frac12\|Av-\xi\|_2^2+\frac12\|v_t\|_2^2
       \quad \text{for } t\in(0,T_{\max}).
  \end{equation}
\end{lemma}
\end{samepage}

\begin{proof}
  Since $A\xi=u$, self-adjointness of $A$ gives
  \begin{equation*}
    \int_\Omega uv=\langle A\xi,v\rangle
                   =\langle\xi,Av\rangle.
  \end{equation*}
  Completing the square therefore yields
  \begin{equation*}
    \frac12\|Av\|_2^2-\int_\Omega uv
      =\frac12\|Av-\xi\|_2^2-\frac12\|\xi\|_2^2.
  \end{equation*}
  Substituting this identity into
  \eqref{sym: energy and dissipation rate} proves
  \eqref{eq:4d-energy-square}.
\end{proof}

We estimate $\|\xi\|_2^2$ in \eqref{eq:4d-energy-square} by a radial resolvent bound.

\begin{lemma}
  \label{lem:4d-resolvent-energy}
  There exists $C_R>0$ such that for any
  $(u_0,v_0,w_0)\in\mathcal I(m,R)$,
  \begin{equation}
    \label{eq:4d-reduced-energy}
    \mathcal F\ge\int_\Omega u\ln u
      -\frac1{16\pi^2}\int_0^R\frac{M(r)^2}{r}\dd r-C_Rm^2
      \quad\text{for }t\in(0,T_{\max}).
  \end{equation}
\end{lemma}

\begin{proof}
  Let $M_{\xi}(r)=\int_{B_r}\xi$. Since $-\Delta\xi+\xi=u$,
  the radial equation is $(r^3\xi_r)'=r^3(\xi-u)$.
  Radial smoothness gives $\lim_{r\searrow0}r^3\xi_r(r)=0$, so
  integration from zero to $r$ yields
  \begin{equation*}
    r^3\xi_r(r)=\int_0^r s^3(\xi(s)-u(s))\dd s
      =\frac{M_{\xi}(r)-M(r)}{2\pi^2}
      \quad\text{for }r\in(0,R].
  \end{equation*}
  Integrating $\xi_r$ from $r$ to $R$, we obtain
  \begin{align*}
    \xi(r)&=\xi(R)+\frac1{2\pi^2}\int_r^R
           \frac{M(s)-M_{\xi}(s)}{s^3}\dd s
       \le \xi(R)+P_0(r),\\
    P_0(r)&\coloneq\frac1{2\pi^2}\int_r^R\frac{M(s)}{s^3}\dd s
    \quad\text{for }r\in(0,R].
\end{align*}
  Fubini's theorem gives
  \begin{equation*}
    \int_\Omega P_0=\frac14\int_0^R sM(s)\dd s
        \le\frac{mR^2}{8}.
  \end{equation*}
  Integration by parts and Fubini's theorem give
  \begin{align*}
    \|P_0\|_2^2
    &=2\pi^2\int_0^Rr^3P_0(r)^2\dd r
      = - 6\pi^2\int_0^Rr^3P_0(r)^2\dd r
      - 4\pi^2\int_0^Rr^4P_0'(r)P_0(r)\dd r\\
      &= - \pi^2\int_0^Rr^4P_0'(r)P_0(r)\dd r
    =\frac1{4\pi^2}\int_0^R rM(r)
        \int_r^R \frac{M(s)}{s^3}\dd s\dd r\\
    &=\frac1{4\pi^2}\int_0^R
        \frac{M(s)}{s^3}\int_0^s rM(r)\dd r\dd s
        \le\frac1{8\pi^2}\int_0^R\frac{M(s)^2}{s}\dd s.
\end{align*}
  The last inequality uses $M(r)\le M(s)$ for $r\le s$.
  The Neumann estimate in \cite[Lemma~23]{Brezis1973} and density
  give the extension $A^{-1}\colon L^1(\Omega)\to W^{1,1}(\Omega)$ with
  \begin{equation}
    \label{eq:neumann-resolvent-W11}
    \|A^{-1}\zeta\|_{W^{1,1}(\Omega)}\le C_R\|\zeta\|_1
    \quad\text{for all }\zeta\in L^1(\Omega).
  \end{equation}
  Since $\xi$ is nonnegative and radial, its trace is constant.
  The trace theorem \cite[Theorem~5.36]{Adams2003}
  (see also \cite[Section~5.5, Theorem~1]{Evans2010}) and
  \eqref{eq:neumann-resolvent-W11} with $\zeta=u$ give
  \begin{equation}
    \label{eq:4d-xi-boundary}
    2\pi^2R^3\xi(R)=\|\xi\|_{L^1(\partial\Omega)}
      \le C_R\|\xi\|_{W^{1,1}(\Omega)}\le C_Rm.
  \end{equation}
  Using $0\le\xi\le\xi(R)+P_0$ and the preceding estimates, we obtain
  \begin{equation}
    \|\xi\|_2^2
    \le\|P_0\|_2^2+2\xi(R)\int_\Omega P_0+|\Omega|\xi(R)^2
    \le\frac1{8\pi^2}\int_0^R\frac{M(r)^2}{r}\dd r+C_Rm^2.
    \label{eq:4d-resolvent-energy}
  \end{equation}
  Inserting \eqref{eq:4d-resolvent-energy} into \eqref{eq:4d-energy-square}
  and discarding the square terms give \eqref{eq:4d-reduced-energy}.
\end{proof}

The critical entropy inequality absorbs the inner mass contribution
in \eqref{eq:4d-reduced-energy}.

\begin{proof}[Proof of \cref{lem:4d-energy-reduction}]
  Set $\mu=\min\{m,64\pi^2\}$. Choose $r_*\in(0,R]$ so that
  $M(r_*)=\mu$ and set
  $u_{\mathrm{i}}=u\mathds{1}_{B_{r_*}}$ and
  $u_{\mathrm{o}}=u\mathds{1}_{\Omega\setminus B_{r_*}}$.
  The inner mass function is $M_{\mathrm{i}}=\min\{M,64\pi^2\}$.

  Since $u_{\mathrm{i}}$ is bounded, nonnegative, and radial, with
  integral $\mu$, \cref{prop:radial-entropy} with $n=4$ and
  $4/\mu\ge1/(16\pi^2)$ gives
  \begin{equation*}
    \int_\Omega u_{\mathrm{i}}\ln u_{\mathrm{i}}
       -\frac1{16\pi^2}\int_0^R\frac{M_{\mathrm{i}}(r)^2}{r}\dd r
    \ge\mu\left(\ln\frac{\mu}{|\Omega|}-1\right)\ge-|\Omega|.
  \end{equation*}
  The disjoint supports give an exact entropy splitting, and
  \begin{equation*}
    0\le M(r)^2-M_{\mathrm{i}}(r)^2
      \le(m+64\pi^2)(M(r)-64\pi^2)_+
      \quad\text{for }r\in[0,R].
  \end{equation*}
  Substituting these bounds and $u_{\mathrm{o}}\ln u_{\mathrm{o}}\ge-1/e$
  into \eqref{eq:4d-reduced-energy} gives \eqref{eq:4d-energy-excess},
  with the two entropy errors and the resolvent error bounded in total by
  $(1+e^{-1})|\Omega|+C_Rm^2\le C_R(1+m^2)$.
\end{proof}

We now bound $Q_R[u]$ in \eqref{eq:4d-energy-excess} using a pointwise
estimate for the local mass excess.

\begin{lemma}
  \label{lem:4d-mass-excess}
  There exists $C_R>0$ such that for any $(u_0,v_0,w_0)\in\mathcal I(m,R)$,
  \begin{equation}
    \label{eq:4d-mass-pointwise}
    (M(r)-64\pi^2)_+
      \le C_Rmr^{2/5}\left(
        1+\frac{\|g\|_2^2+\|A^{-3/4}\mathcal R\|_2^2}{1+m^2}
      \right)^{1/2}
       \quad\text{for }(r,t)\in(0,R]\times(0,T_{\max}).
  \end{equation}
  Consequently, the weighted integral defined in
  \eqref{eq:4d-excess-moment-definition} satisfies
  \begin{equation}
    \label{eq:4d-excess-moment}
    Q_R[u]
    \le C_Rm\left[1+\ln\left(
      1+\frac{\|g\|_2^2+\|A^{-3/4}\mathcal R\|_2^2}{1+m^2}
    \right)\right]
    \quad \text{for } t\in(0,T_{\max}),
  \end{equation}
  whenever $(u_0,v_0,w_0)\in\mathcal I(m,R)$.
\end{lemma}

The pointwise estimate \eqref{eq:4d-mass-pointwise} is based on a
quadratic radial mass identity.

\begin{lemma}
  \label{lem:4d-quadratic-mass}
  Set $z=A^{-2}u=A^{-1}\xi$, so that $h=v-z$. Define
  \begin{equation*}
    U(r)=\int_0^r s^3u(s)\dd s=\frac{M(r)}{2\pi^2},\qquad
    E(r)=\int_0^r s^3(2\xi(s)-z(s))\dd s
    \quad\text{for }r\in[0,R].
  \end{equation*}
  For an accumulated radial function $F\in C([0,R])$, put
  \begin{equation}
    \label{eq:4d-K-operator}
    K_F(r)\coloneq\frac1{4r^2}\int_0^r sF(s)\dd s
            +\frac{r^2}{4}\int_r^R\frac{F(s)}{s^3}\dd s
    \quad\text{for }r\in(0,R].
  \end{equation}
  Then
  \begin{align}
    \label{eq:4d-quadratic-mass}
    r^4u(r)
    ={}&4U(r)-\frac{U(r)^2}{8}
     -\frac18\left(\int_0^r s^5u(s)\dd s\right)
              \left(\int_r^R su(s)\dd s\right)\notag\\
     &+\left(\frac{U(R)}{8R^2}+\frac{z(R)-\xi(R)}4\right)
           \int_0^r s^5u(s)\dd s
       +\int_0^r U'(s)K_E(s)\dd s\notag\\
     &+\int_0^r U'(s)s h_r(s)\dd s
       +\int_0^r s^4\sqrt{u(s)}g(s)\dd s
    \quad\text{for }r\in(0,R].
  \end{align}
\end{lemma}

\begin{proof}
  Since $u$ is bounded, $s^4u(s)\to0$ as $s\searrow0$.
  Integrating $(s^4u(s))'=4s^3u(s)+s^4u_r(s)$ and using
  $u_r=uv_r+\sqrt u\,g$ and $U'(s)=s^3u(s)$ yields
  \begin{align*}
    r^4u(r)
    &=4U(r)+\int_0^r s^4u_r(s)\dd s\\
    &=4U(r)+\int_0^r U'(s)s v_r(s)\dd s
       +\int_0^r s^4\sqrt{u(s)}g(s)\dd s
    \quad\text{for }r\in(0,R].
  \end{align*}
  We decompose $v_r=z_r+h_r$ and first compute the contribution
  of $z_r$. Set $y=\Delta z=z-\xi$.
  Since $A^2z=u$, expansion gives
  $\Delta y=\Delta^2z=u-(2\xi-z)$. Thus
  $(r^3y_r)'=U'-E'$, and radial integration yields
  \begin{equation*}
    r^3y_r(r)=U(r)-E(r),\qquad
    y(r)=y(R)-\int_r^R\frac{U(s)-E(s)}{s^3}\dd s
    \quad\text{for }r\in(0,R].
  \end{equation*}
  Integrating $(r^3z_r)'=r^3y$ from zero and then integrating
  by parts, we obtain
  \begin{equation*}
    rz_r(r)
    =\frac1{r^2}\int_0^r s^3y(s)\dd s
    =\frac{r^2y(r)}4
      -\frac1{4r^2}\int_0^r s\bigl(U(s)-E(s)\bigr)\dd s.
  \end{equation*}
  Substituting the preceding formula for $y(r)$ gives
  \begin{equation}
    \label{eq:4d-quasistatic-gradient}
    rz_r(r)=\frac{y(R)}4r^2-K_U(r)+K_E(r)
    \quad\text{for }r\in(0,R].
  \end{equation}

  To compute the contribution of $K_U$, define
  \begin{equation*}
    P(r)=\int_0^r s^2U'(s)\dd s,\qquad
    H(r)=\int_r^R\frac{U'(s)}{s^2}\dd s-\frac{U(R)}{R^2}.
  \end{equation*}
  Integration by parts in the two integrals defining $K_U$ gives
  \begin{equation*}
    K_U(r)=\frac{U(r)}4-\frac{P(r)}{8r^2}+\frac{r^2H(r)}8.
  \end{equation*}
  Since $P'=r^2U'$ and $H'=-U'/r^2$, it follows that
  \begin{equation*}
    U'K_U=\frac{UU'}4+\frac{P'H+PH'}8
      =\frac18\bigl(U^2+PH\bigr)'.
  \end{equation*}
Integration from zero to $r$ therefore yields
  \begin{equation*}
    \int_0^r U'(s)K_U(s)\dd s
      =\frac{U(r)^2}{8}
       +\frac18\left(\int_0^r s^5u(s)\dd s\right)
                \left(\int_r^R su(s)\dd s-\frac{U(R)}{R^2}\right).
  \end{equation*}
  The boundary term at zero vanishes, thanks to boundedness of $u$.
  Substituting this identity and \eqref{eq:4d-quasistatic-gradient}
  into the starting mass identity, and using $v_r=z_r+h_r$ and
  $y(R)=z(R)-\xi(R)$, proves \eqref{eq:4d-quadratic-mass}.
\end{proof}

The quadratic term $4U-U^2/8$ in \eqref{eq:4d-quadratic-mass}
changes sign as $M$ crosses $64\pi^2$. Bounds on the remaining
terms give \eqref{eq:4d-mass-pointwise}.

\begin{proof}[Proof of \cref{lem:4d-mass-excess}]
  Integrating $A\xi=u$ and using the Neumann condition give
  $\|\xi\|_1=m$ by nonnegativity. Applying
  \cite[Lemma~23]{Brezis1973} to
  the homogeneous Neumann problems $A\xi=u$ and $Az=\xi$ gives
  \begin{equation*}
    \|\xi\|_{10/9}+\|z\|_{10/9}
      \le C_R(\|u\|_1+\|\xi\|_1)\le C_Rm.
  \end{equation*}
  H\"older's inequality therefore yields
  \begin{equation*}
    |E(r)|\le\frac1{2\pi^2}\|2\xi-z\|_{10/9}|B_r|^{1/10}
      \le C_Rmr^{2/5}
      \quad\text{for }r\in(0,R].
  \end{equation*}
  Integration of the two kernels in
  \eqref{eq:4d-K-operator} therefore gives
  \begin{equation}
    \label{eq:4d-lower-order-kernel}
    |K_E(r)|\le C_Rmr^{2/5}
    \quad\text{for }r\in(0,R].
  \end{equation}
  For the residual potential, set $\eta=A^{-3/4}\mathcal R$.
  Since $h=A^{-2}\mathcal R=A^{-5/4}\eta$,
  \begin{equation*}
    \Delta h=A^{-5/4}\eta-A^{-1/4}\eta.
  \end{equation*}
  The semigroup representation \eqref{eq:4d-negative-powers} and
  the $L^2$--$L^{5/2}$ heat estimate
  \cite[Lemma~1.3(i)]{Winkler2010b}, after separating the spatial
  mean, give
  \begin{align*}
    \|A^{-1/4}\psi\|_{5/2}+\|A^{-5/4}\psi\|_{5/2}
    &\le C_R\|\psi\|_2\int_0^\infty
      (\tau^{-3/4}+\tau^{1/4})(1+\tau^{-1/5})e^{-\tau}\dd\tau\\
    &\le C_R\|\psi\|_2
    \quad\text{for }\psi\in L^2(\Omega).
  \end{align*}
  Applying this bound to $\psi=\eta$ gives
  $\|\Delta h\|_{5/2}\le C_R\|\eta\|_2$.
  Radial integration and H\"older's inequality now imply
  \begin{equation}
    \label{eq:4d-residual-gradient}
    |rh_r(r)|
      =\left|\frac1{r^2}\int_0^r s^3\Delta h(s)\dd s\right|
      \le C_Rr^{2/5}\|\eta\|_2
      \quad\text{for }r\in(0,R].
  \end{equation}

  Since $z=A^{-1}\xi$ and $\|\xi\|_1=m$, the trace calculation
  in \eqref{eq:4d-xi-boundary} also gives $z(R)\le C_Rm$. Hence
  \begin{equation}
    \label{eq:4d-quasistatic-boundary}
    \xi(R)+z(R)\le C_Rm.
  \end{equation}
  In \eqref{eq:4d-quadratic-mass}, the exterior mass term is
  nonpositive. The boundary coefficient is bounded by
  \eqref{eq:4d-quasistatic-boundary}, and
  $\int_0^r s^5u\dd s\le r^2U(r)$. For the two integrals against
  $U'$, use \eqref{eq:4d-lower-order-kernel} and
  \eqref{eq:4d-residual-gradient} with $\int_0^r U'(s)\dd s=U(r)$.
  Cauchy--Schwarz bounds the flux term by
  \begin{equation*}
    \left|2\pi^2\int_0^r s^4\sqrt u\,g\dd s\right|
      =\left|\int_{B_r}|x|\sqrt u\,g\right|
      \le r\sqrt{M(r)}\,\|g\|_2
    \quad\text{for }r\in(0,R].
  \end{equation*}
  Multiplying \eqref{eq:4d-quadratic-mass} by $2\pi^2$
  and using $r,r^2\le C_Rr^{2/5}$ gives
  \begin{align*}
    0\le rM'(r)
    &\le 4M(r)\left(1-\frac{M(r)}{64\pi^2}\right)\\
    &\quad+C_Rr^{2/5}
      \bigl[M(r)(m+\|\eta\|_2)+\sqrt{M(r)}\,\|g\|_2\bigr]
    \quad\text{for }r\in(0,R].
  \end{align*}
  When $M(r)>64\pi^2$, divide by $M(r)$ to obtain
  \begin{equation*}
    M(r)-64\pi^2
    \le C_Rr^{2/5}\left(m+\|\eta\|_2+\frac{\|g\|_2}{\sqrt{M(r)}}\right)
    \quad\text{for }r\in(0,R]\text{ with }M(r)>64\pi^2.
  \end{equation*}
  In this range $m\ge M(r)>64\pi^2>1$, so
  $1+m^2\le2m^2$. By Cauchy--Schwarz the quantity in parentheses
  is at most $\sqrt{6}\,m\bigl(1+(\|g\|_2^2+\|\eta\|_2^2)/(1+m^2)\bigr)^{1/2}$.
  For $M(r)\le64\pi^2$ the positive part vanishes.
  This proves \eqref{eq:4d-mass-pointwise}.

  Set $\Lambda=C_R\bigl(1+(\|g\|_2^2+\|\eta\|_2^2)/(1+m^2)\bigr)^{1/2}\ge1$,
  with $C_R$ chosen as in \eqref{eq:4d-mass-pointwise}, and let
  $r_0=\min\{R,\Lambda^{-5/2}\}$. Using
  $(M(r)-64\pi^2)_+\le m\min\{\Lambda r^{2/5},1\}$ gives
  \begin{equation*}
    Q_R[u]
    \le m\left(\Lambda\int_0^{r_0}r^{-3/5}\dd r
                 +\ln\frac R{r_0}\right)
    \le m\left[\frac52+\max\{0,\ln(R\Lambda^{5/2})\}\right].
  \end{equation*}
  Substituting $\Lambda$ gives \eqref{eq:4d-excess-moment}.
\end{proof}

\begin{proof}[Proof of \cref{prop:4d-static-log-energy}]
  Insert \eqref{eq:4d-excess-moment} into \eqref{eq:4d-energy-excess}.
  Since $m(m+64\pi^2)\le(1+64\pi^2)(1+m^2)$,
  \eqref{eq:4d-static-log-energy} follows upon enlarging $C_R$ and
  taking the positive part.
\end{proof}

\subsection{Temporal estimates for the elliptic residual}
\label{sec:4d-closed-proof}

 The main result of this subsection is
the following time-integrated estimate for the elliptic residual in \eqref{eq:4d-static-log-energy}.

\begin{lemma}
  \label{lem:4d-residual-time}
  Under the hypotheses of \cref{prop:4d-closed-spacetime}, there
  exists $C_R>0$ such that
  \begin{equation}
    \int_0^s \|A^{-3/4}\mathcal R(t)\|_2^2\dd t
        \le C_R\left[(m+L)^2+(1+\sqrt m)
                   \int_0^s\mathcal D(t)\dd t\right]
       \quad\text{for }s\in(0,T_{\max}).
    \label{eq:4d-residual-time}
  \end{equation}
\end{lemma}

We prepare uniform estimates of signals.

\begin{lemma}
  \label{lem:4d-weak-signal-bounds}
  Under the assumptions of \cref{prop:4d-closed-spacetime},
  there exists $C_R>0$ such that
  \begin{equation}
    \label{eq:4d-endpoint-signal-bounds}
    \|w(t)\|_{W^{1,5/4}}
      +\|v(t)\|_{W^{2,7/4}}
      +\|Av(t)\|_{7/4}+\|v_t(t)\|_{7/4}\le C_R(m+L)
    \quad\text{for }t\in(0,T_{\max}).
  \end{equation}
\end{lemma}

\begin{proof}
  The constant $C_R$ may vary from line to line.
  The Neumann heat semigroup is uniformly bounded on
  $W^{1,5/4}(\Omega)$. Indeed,
  $D_{5/4}(A^{1/2})=W^{1,5/4}(\Omega)$ with equivalent norms
  by \cite[Theorem~2.2]{Grubb2016}. By
  $L^{5/4}$-contractivity and commutation of the fractional power
  with the semigroup \cite[Theorem~2.6.13(b)]{Pazy1983},
  \begin{equation*}
    \|e^{-\tau A}\phi\|_{W^{1,5/4}}
      \le C_Re^{-\tau}\|\phi\|_{W^{1,5/4}}
      \quad\text{for }\tau>0,\quad
      \phi\in W^{1,5/4}(\Omega).
  \end{equation*}
  The $L^1$ smoothing and gradient estimates in
  \cite[Lemma~1.3(i)--(ii)]{Winkler2010b} also give
  \begin{equation*}
    \|e^{-\tau A}\phi\|_{W^{1,5/4}}
      \le C_R(1+\tau^{-9/10})e^{-\tau}\|\phi\|_1
      \quad\text{for }\tau>0,\quad\phi\in L^1(\Omega).
  \end{equation*}
  The variation-of-constants formula reads
  \begin{equation*}
    w(t)=e^{-tA}w_0
      +\int_0^t e^{-(t-\sigma)A}u(\sigma)\dd\sigma
    \quad\text{for }t\in(0,T_{\max}).
  \end{equation*}
   Applying the semigroup estimates, mass preservation $\|u\|_1=m$ and
  \eqref{eq:4d-initial-weak-bound} then gives
  \begin{equation}
    \label{eq:4d-w-weak-bound}
    \|w(t)\|_{W^{1,5/4}}
    \le C_R L+C_Rm\int_0^t
        (1+\tau^{-9/10})e^{-\tau}\dd\tau
    \le C_R(m+L)
    \quad\text{for }t\in(0,T_{\max}).
  \end{equation}

  To estimate the Hessian of $v$, fix $i\in\{1,2,3,4\}$ and
  differentiate its equation in $x_i$. Radial symmetry and
  $\partial_\nu v=0$ imply $v_{x_i}=0$ on $\partial B_R$.
  Thus
  \begin{equation*}
    \begin{cases}
      \partial_t v_{x_i}=\Delta v_{x_i}-v_{x_i}+w_{x_i}
        &\text{for }x\in\Omega,\ t\in(0,T_{\max}),\\
      v_{x_i}=0
        &\text{for }x\in\partial\Omega,\ t\in(0,T_{\max}),\\
      v_{x_i}(x,0)=v_{0x_i}(x)
        &\text{for }x\in\Omega.
    \end{cases}
  \end{equation*}
  Let $\Delta_D$ denote the Dirichlet Laplacian on $B_R$ and
  set $S_D(\tau)=e^{\tau(\Delta_D-1)}$ for $\tau\ge0$.
  By \cite[Lemma~2.4(i)]{Mu2024},
  \begin{equation*}
    \|\nabla S_D(\tau)\phi\|_{7/4}
      \le C_R(1+\tau^{-67/70})e^{-\tau}\|\phi\|_{5/4}
      \quad\text{for }\tau>0,\quad\phi\in L^{5/4}(\Omega).
  \end{equation*}
  Moreover, \cite[Theorem~2.2 and Corollary~2.3]{Grubb2016} gives
  $D_{7/4}((1-\Delta_D)^{1/2})=W^{1,7/4}_0(B_R)$ with equivalent
  norms. Commutation of this fractional power with $S_D(\tau)$,
  as in \cite[Theorem~2.6.13(b)]{Pazy1983},
  and $L^{7/4}$-contractivity of $e^{\tau\Delta_D}$ imply
  \begin{equation*}
    \|S_D(\tau)\phi\|_{W^{1,7/4}}
      \le C_Re^{-\tau}\|\phi\|_{W^{1,7/4}}
      \quad\text{for }\tau>0,\quad
      \phi\in W^{1,7/4}_0(\Omega).
  \end{equation*}
  The variation-of-constants formula is
  \begin{equation*}
    v_{x_i}(t)=S_D(t)v_{0x_i}
      +\int_0^t S_D(t-\sigma)w_{x_i}(\sigma)\dd\sigma
    \quad\text{for }t\in(0,T_{\max}).
  \end{equation*}
  Here $v_{0x_i}\in W^{1,7/4}_0(\Omega)$, since radial symmetry
  and the Neumann condition imply $v_{0x_i}=0$ on $\partial\Omega$.
  Applying the two Dirichlet semigroup estimates gives
  \begin{align*}
      \|\nabla v_{x_i}(t)\|_{7/4}
      &\le C_R\|v_0\|_{W^{2,7/4}}
         +C_R\sup_{0<\sigma<t}\|\nabla w(\sigma)\|_{5/4}
            \int_0^t(1+\tau^{-67/70})e^{-\tau}\dd\tau\\
      &\le C_R(m+L)
      \quad\text{for }t\in(0,T_{\max}),\quad i\in\{1,2,3,4\}.
\end{align*}
  Integrating the two signal equations gives
  \begin{equation*}
    \|v(t)\|_1
       \le\max\left\{m,\int_\Omega v_0,\int_\Omega w_0\right\}
       \le C_R(m+L)
       \quad\text{for }t\in(0,T_{\max}),
  \end{equation*}
by comparison,
  nonnegativity, and \eqref{eq:4d-initial-weak-bound},
  as done in \cite[Lemma~2.3]{Mao2025b}.
  By interpolation, $\|v(t)\|_{W^{2,7/4}}
  \le C_R(\|D^2v(t)\|_{7/4}+\|v(t)\|_1)\le C_R(m+L)$
  for $t\in(0,T_{\max})$. Since $A=1-\Delta_N$ and $v_t=w-Av$,
  $\|Av(t)\|_{7/4}+\|v_t(t)\|_{7/4}
  \le C_R\|v(t)\|_{W^{2,7/4}}+\|w(t)\|_{7/4}$.
  Together with \eqref{eq:4d-w-weak-bound},
  and Sobolev embedding $W^{1,5/4}(\Omega) \hookrightarrow L^{7/4}(\Omega)$,
  this proves \eqref{eq:4d-endpoint-signal-bounds}.
\end{proof}

\begin{proof}[Proof of \cref{lem:4d-residual-time}]
  Write the cross functional as
  \begin{equation}
    \label{eq:4d-cross-functional}
    \mathcal J(t)
      =\langle v_t,A^{-3/2}\mathcal R\rangle
      =\langle v_t,A^{1/2}v-A^{-3/2}u\rangle
      \quad\text{for }t\in(0,T_{\max}).
  \end{equation}
  For the endpoint bound, apply the Neumann heat estimate
  \cite[Lemma~1.3(i)]{Winkler2010b} with target space $L^{7/3}$,
  separating the spatial mean. Inserting it into
  \eqref{eq:4d-negative-powers} and using $\|u\|_1=m$ gives
  \begin{align}
    \|A^{-3/2}u\|_{7/3}
    &\le C_Rm\int_0^\infty
       \tau^{1/2}(1+\tau^{-8/7})e^{-\tau}\dd\tau
    \le C_Rm,\label{eq:4d-fractional-L1}\\
    \|A^{-1/2}Av\|_{7/3}
    &\le C_R\|Av\|_{7/4}\int_0^\infty
       \tau^{-1/2}(1+\tau^{-2/7})e^{-\tau}\dd\tau
    \le C_R\|Av\|_{7/4}
    \quad\text{for }t\in(0,T_{\max}).\label{eq:4d-fractional-Av}
  \end{align}
  By H\"older's inequality and
  \eqref{eq:4d-endpoint-signal-bounds},
  \begin{equation}
    \label{eq:4d-cross-endpoint}
    |\mathcal J(t)|
      \le \|v_t\|_{7/4}
          \bigl(\|A^{-1/2}Av\|_{7/3}
                +\|A^{-3/2}u\|_{7/3}\bigr)
      \le C_R(m+L)^2
      \quad\text{for }t\in(0,T_{\max}).
  \end{equation}
  Differentiating the second expression in
  \eqref{eq:4d-cross-functional} gives
  \begin{align*}
    \mathcal J'
    &=\langle v_{tt},A^{1/2}v-A^{-3/2}u\rangle\\
    &\quad+\langle v_t,A^{1/2}v_t\rangle
       -\langle A^{-3/2}v_t,u_t\rangle\\
    &=-\|A^{-3/4}\mathcal R\|_2^2
       -2\langle A^{1/4}v_t,A^{-3/4}\mathcal R\rangle\\
    &\quad+\|A^{1/4}v_t\|_2^2
       -\langle A^{-3/2}v_t,u_t\rangle
       \quad\text{for }t\in(0,T_{\max}).
  \end{align*}
  Here we used $A^{1/2}v-A^{-3/2}u=A^{-3/2}\mathcal R$,
  the identity $v_{tt}=-2Av_t-\mathcal R$ from
  \eqref{eq:4d-damped-signal}, the self-adjointness of the powers of $A$,
  and their composition rule \cite[Theorem~2.6.8(d)]{Pazy1983}.
  Since $A^{-3/2}v_t$ is radial, the first equation in
  \eqref{eq:4d-residual-evolution} gives, by radial integration by parts,
  \begin{equation*}
    -\langle A^{-3/2}v_t,u_t\rangle
    =\int_\Omega\sqrt u\,g\,\partial_r(A^{-3/2}v_t)
    \quad\text{for }t\in(0,T_{\max}).
  \end{equation*}
  The endpoint terms vanish by the zero flux at $r=R$ and radial
  regularity at $r=0$.
  To bound the radial derivative, write
  $A^{-3/2}v_t=A^{-2}(A^{1/2}v_t)$ and apply the
  $L^2$--$L^\infty$ gradient estimate
  \cite[Lemma~1.3(ii)]{Winkler2010b} in \eqref{eq:4d-negative-powers}.
  This yields
  \begin{align}
    \|\nabla A^{-3/2}v_t\|_\infty
    &\le C_R\|A^{1/2}v_t\|_2
       \int_0^\infty\tau(1+\tau^{-3/2})e^{-\tau}\dd\tau
       \notag\\
    &\le C_R\|A^{1/2}v_t\|_2
       \quad\text{for }t\in(0,T_{\max}).
    \label{eq:4d-fractional-gradient}
  \end{align}
  Cauchy--Schwarz, mass conservation, and
  \eqref{eq:4d-dissipation-operator} now give
  \begin{equation*}
    \left|\int_\Omega\sqrt u\,g\,
                     \partial_r(A^{-3/2}v_t)\right|
      \le C_R\sqrt m\,\|g\|_2\|A^{1/2}v_t\|_2
      \le C_R\sqrt m\,\mathcal D
      \quad\text{for }t\in(0,T_{\max}).
  \end{equation*}
  Moreover,
  \begin{equation*}
    2|\langle A^{1/4}v_t,A^{-3/4}\mathcal R\rangle|
       \le\frac12\|A^{-3/4}\mathcal R\|_2^2+2\|A^{1/4}v_t\|_2^2,
    \qquad
    \|A^{1/4}v_t\|_2\le\|A^{1/2}v_t\|_2,
  \end{equation*}
  because $A\ge1$. Therefore
  \begin{equation}
    \label{eq:4d-residual-differential}
    \mathcal J'(t)+\frac12\|A^{-3/4}\mathcal R(t)\|_2^2
       \le C_R(1+\sqrt m)\mathcal D(t)
       \quad\text{for }t\in(0,T_{\max}).
  \end{equation}
  Integrating \eqref{eq:4d-residual-differential} from $\sigma$ to $s$
  and using \eqref{eq:4d-cross-endpoint} at both endpoints yields
  \begin{align*}
    \frac12\int_{\sigma}^s\|A^{-3/4}\mathcal R(t)\|_2^2\dd t
       &\le \mathcal J(\sigma)-\mathcal J(s)
            +C_R(1+\sqrt m)\int_{\sigma}^s\mathcal D(t)\dd t\\
       &\le C_R\left[(m+L)^2+(1+\sqrt m)
                 \int_0^s\mathcal D(t)\dd t\right]
         \quad\text{for }0<\sigma<s<T_{\max}.
\end{align*}
  By \eqref{eq:4d-integrated-energy}, the dissipation integral on
  the right is finite.
  \eqref{eq:4d-residual-time} follows from letting $\sigma\searrow0$ by monotone convergence.
\end{proof}

\subsection{The closed spacetime estimate and its consequence}
\label{sec:4d-blowup-proof}

We combine the spatial estimate of \cref{sec:4d-spatial-proof}
with the temporal residual control of \cref{sec:4d-closed-proof}
to prove \cref{prop:4d-closed-spacetime}. The resulting closed
spacetime estimate yields the low-energy blow-up criterion by an
auxiliary-function argument. Finally, the Fujie--Senba construction
provides initial data satisfying this criterion at every
supercritical mass.

\begin{proof}[Proof of \cref{prop:4d-closed-spacetime}]
  Set $H(t)=(\|g(t)\|_2^2+\|A^{-3/4}\mathcal R(t)\|_2^2)/(1+m^2)$
  for $t\in(0,T_{\max})$. By \eqref{eq:4d-dissipation-operator}
  and \eqref{eq:4d-residual-time},
  \begin{align*}
    \int_0^s H(t)\dd t
    &\le\frac{C_R}{1+m^2}\left[(m+L)^2
        +(1+\sqrt m)\int_0^s\mathcal D(t)\dd t\right]\\
    &\le C_R\left[1+\frac{L^2}{1+m^2}
        +\int_0^s\mathcal D(t)\dd t\right]\\
    &\le\lambda\left(1+\int_0^s\mathcal D(t)\dd t\right)<\infty
    \quad\text{for }s\in(0,T_{\max}),
  \end{align*}
  where $\lambda=C_R(1+L^2/(1+m^2))\ge1$ after enlarging $C_R$.
  The second inequality uses $(m+L)^2\le2m^2+2L^2$ and
  $1+\sqrt m\le2(1+m^2)$.
  By \eqref{eq:4d-static-log-energy} and Jensen's inequality with
  the probability measure $\dd t/s$ on $(0,s)$,
  \begin{align*}
    \int_0^s(-\mathcal F(t))_+\dd t
    &\le C_R(1+m^2)s\left[1+\ln\left(
          1+\frac{1}{s}\int_0^s H(t)\dd t\right)\right]\\
    &\le C_R(1+m^2)s\left[1+\ln\left(
          1+\lambda\frac{1+\int_0^s\mathcal D(t)\dd t}{s}\right)\right]
    \quad\text{for }s\in(0,T_{\max}).
  \end{align*}
  Since $1+\lambda x\le\lambda(1+x)$ for $\lambda\ge1$,
  \begin{equation*}
    1+\ln(1+\lambda x)
      \le(1+\ln\lambda)\bigl[1+\ln(1+x)\bigr]
      \quad\text{for }x\ge0.
  \end{equation*}
  The definition of $\lambda$ gives
  \begin{equation*}
    1+\ln\lambda\le C_R\ln\left(e+\frac{L^2}{1+m^2}\right).
  \end{equation*}
  Applying these bounds to the preceding estimate proves
  \eqref{eq:4d-intro-closed}.
 \end{proof}

\begin{proof}[Proof of \cref{thm:4d-low-energy-blowup}~(i)]
  Choose $C=C(m,R,L)\ge1$ such that \eqref{eq:4d-intro-closed} holds under assumptions there,
  and set
  \begin{equation*}
    \ell\coloneq2C\ln(4C),\qquad
    K\coloneq C+\ell,\qquad T\coloneq2.
  \end{equation*}
  Assume \eqref{eq:4d-low-energy-data} and suppose for contradiction
  that $T_{\max}>T$.
  Following the auxiliary-function construction in
  \cite{Mao2025b}, define
  \begin{equation*}
    \Psi(s)\coloneq
       \int_0^s\bigl[-\mathcal F(t)+\mathcal F(0)+\ell\bigr]\dd t
       \quad\text{for }s\in[0,T_{\max}).
  \end{equation*}
  By \eqref{eq:4d-integrated-energy},
  \begin{equation*}
    \Psi'(s)=\mathcal F(0)-\mathcal F(s)+\ell
       =\int_0^s\mathcal D(t)\dd t+\ell\ge\ell,
    \qquad
    \Psi(s)\ge\ell s
    \quad\text{for }s\in(0,T_{\max}).
  \end{equation*}
  Since $\mathcal F(t)\le\mathcal F(0)<0$, the closed estimate
  \eqref{eq:4d-intro-closed} yields
  \begin{align*}
    \Psi(s)
    &\le Cs\left[1+\ln\left(
           1+\frac{\Psi'(s)-\ell+1}{s}\right)\right]
         +\bigl(\mathcal F(0)+\ell\bigr)s\\
    &\le Cs\ln\left(1+\frac{\Psi'(s)}s\right)
       \quad\text{for }s\in(0,T_{\max}),
  \end{align*}
  where the last step uses $\ell\ge1$ and
  $\mathcal F(0)+\ell+C\le0$.
  In particular, on $[1,2]$,
  \begin{equation*}
    \Psi(s)\le2C\ln\bigl(1+\Psi'(s)\bigr).
  \end{equation*}
  Inverting this inequality gives
  \begin{equation*}
    \Psi'(s)\ge e^{\Psi(s)/(2C)}-1
       \ge\frac12e^{\Psi(s)/(2C)}
       \quad\text{for }s\in[1,2],
  \end{equation*}
  because $\Psi(s)\ge\ell$ and $e^{\ell/(2C)}=4C\ge4$.
  Consequently,
  \begin{equation*}
    \frac{\dd}{\dd s}e^{-\Psi(s)/(2C)}
       =-\frac{\Psi'(s)}{2C}e^{-\Psi(s)/(2C)}
       \le-\frac1{4C}
       \quad\text{for }s\in[1,2].
  \end{equation*}
  Integration from $1$ to $2$ gives the contradiction
  \begin{equation*}
    0<e^{-\Psi(2)/(2C)}
       \le e^{-\Psi(1)/(2C)}-\frac1{4C}
       \le e^{-\ell/(2C)}-\frac1{4C}=0.
  \end{equation*}
  Hence $T_{\max}\le2=T$.

  \begin{samepage}
  \emph{Proof of \cref{thm:4d-low-energy-blowup}~(ii).}
  Fix $m>64\pi^2$. We use the construction of Fujie--Senba
  \cite[Section~4]{Fujie2019}, with
  $\chi=\tau_1=\tau_2=1$ in their notation.
  Fix $0<r_1<r_2<\min\{R/2,1\}$ and a radial cutoff
  $\phi\in C_c^\infty(B_{r_2})$ with $0\le\phi\le1$,
  $\phi=1$ on $B_{r_1}$, and $x\cdot\nabla\phi(x)\le0$.
  Their profiles are
  \begin{align*}
    \bar u_\lambda(x)
      &=\frac{384\lambda^4}{(1+\lambda^2|x|^2)^4},\\*
    \bar v_\lambda(x)
      &=4\ln\frac{1+\lambda^2r_2^2}{1+\lambda^2|x|^2}
        +\ln384,\\*
    \bar w_\lambda(x)
      &=\frac{16\lambda^2(2+\lambda^2|x|^2)}
              {(1+\lambda^2|x|^2)^2}
      \quad\text{for }x\in\mathbb R^4,\quad\lambda\ge1.
  \end{align*}
  \end{samepage}
  Define
  \begin{equation*}
    a_\lambda =\frac{m}{\int_\Omega\bar u_\lambda\phi},\qquad
    (u_{0,\lambda},v_{0,\lambda},w_{0,\lambda})(x)
       =a_\lambda\phi(x)
        (\bar u_\lambda,\bar v_\lambda,\bar w_\lambda)(x)
      \quad\text{for }x\in\Omega,\quad\lambda\ge1.
  \end{equation*}
  Since $\bar v_\lambda\ge\ln384>0$ on $B_{r_2}$,
  these data are nonnegative, smooth, and radial, have cell mass
  $m$, and vanish near $\partial\Omega$. In particular,
  $(u_{0,\lambda},v_{0,\lambda},w_{0,\lambda})\in\mathcal I(m,R)$.
  The identity
  $\int_{\mathbb R^4}\bar u_\lambda=64\pi^2$ and the change of
  variables $y=\lambda x$ imply
  \begin{equation*}
    a_\lambda>\frac{m}{64\pi^2}>1,\qquad
    a_\lambda\to\frac{m}{64\pi^2}
    \quad\text{as }\lambda\to\infty.
  \end{equation*}
  Thus $a_\lambda$ is uniformly bounded for $\lambda\ge1$.

  The explicit profiles and the fixed cutoff give
  \begin{equation*}
    |\nabla w_{0,\lambda}(x)| \le C|x|^{-3},\qquad
    |D^2v_{0,\lambda}(x)| \le C|x|^{-2}
    \quad\text{for }x\in\Omega\setminus\{0\},\quad\lambda\ge1.
  \end{equation*}
  Since $3\cdot(5/4)<4$ and $2\cdot(7/4)<4$, radial integration
  and Poincar\'e inequality yield a constant $L_*>0$,
  independent of $\lambda$, such that
  \begin{equation*}
    \|w_{0,\lambda}\|_{W^{1,5/4}(\Omega)}
      +\|v_{0,\lambda}\|_{W^{2,7/4}(\Omega)}\le L_*
    \quad\text{for }\lambda\ge1.
  \end{equation*}
  The energy calculation in \cite[(4.6)--(4.10)]{Fujie2019}
  includes the kinetic term in \eqref{eq:4d-initial-energy}
  and gives, for these data,
  \begin{align*}
    \mathcal F(0)
      &\le-256\pi^2a_\lambda(a_\lambda-1)\ln\lambda+C\\
      &\le-4m\left(\frac{m}{64\pi^2}-1\right)\ln\lambda+C
        \to-\infty
      \quad\text{as }\lambda\to\infty.
  \end{align*}
  Choose $\lambda$ large enough that the energy is at most
  $-K(m,R,L_*)$. Part~(i), applied with $L=L_*$, then gives
  finite-time blow-up. Taking a positive time before blow-up
  as the new initial time gives smooth, strictly positive radial
  data of the same mass whose solution still blows up in finite time.
\end{proof}

\appendix

\section{Optimality of the logarithmic bound in two dimensions}
\label{sec:2d-supercritical-data}

We prove the optimality of the logarithmic growth in
\eqref{eq:2d-sobolev-energy-dissipation} by pairing a concentrating
density with a signal obtained from a Neumann--Poisson equation.
For an earlier low-energy construction, see Mizoguchi--Winkler
\cite[Lemma~8.2]{Mizoguchi2014}.

\begin{proposition}
\label{prop:2d-concentrating-sequence}
Fix $m>8\pi$ and $R>0$. Let $\Omega=B_R\subset\mathbb R^2$.
There exists a sequence $(u_k,v_k)_{k\ge1}$ of smooth radial function pairs on $B_R$ satisfying
\begin{equation}
\label{eq:2d-family-mass-norm}
(u_k,v_k)\in\mathcal S(m,R)
\quad\text{for }k\ge1,
\end{equation}
and
\begin{equation}
\label{eq:2d-family-energy}
\mathcal F_{\mathrm{KS}}(u_k,v_k)
=2m\left(1-\frac{m}{8\pi}\right)\ln k+O(1)
\to-\infty\quad\text{as }k\to\infty.
\end{equation}
Moreover,
\begin{equation}
\label{eq:2d-family-log-dissipation}
\ln\bigl(1+\mathcal D_{\mathrm{KS}}(u_k,v_k)\bigr)
=\frac{4(m-4\pi)}{m-8\pi}\ln k+O(1)
\quad\text{as }k\to\infty.
\end{equation}
Consequently, the function $\ln(1+s)$ in
\eqref{eq:2d-sobolev-energy-dissipation} cannot be replaced by
any nonnegative function $h$ satisfying $h(s)=o(\ln(1+s))$
as $s\to\infty$.
\end{proposition}

\begin{proof}
Fix $\rho\in(0,R)$ and set $b=m/(m-8\pi)$ and
$a=1-(1+b)^{-1/2}\in(0,1)$.
Choose a nonnegative radial $\psi\in C_c^\infty(B_1)$ with
$\|\psi\|_{L^2(B_1)}=1$, put $\phi=\psi^2$, and define
$\phi_\sigma(x)=\sigma^{-2}\phi(x/\sigma)$ for $x\in\mathbb R^2$ and
$\sigma>0$. For $k\ge1$ let
$\varepsilon_k=\rho k^{-(1+b)}$ and $\delta_k=1/k$. Define
\begin{equation}
\label{eq:2d-concentrating-data}
u_k(x)=m(1-\delta_k)\phi_{\varepsilon_k}(x)+\frac{m\delta_k}{|\Omega|}
\quad\text{for }x\in\Omega,\ k\ge1.
\end{equation}
Let $v_k$ be the radial solution, normalized by $v_k(R)=1$, of
\begin{equation}
\label{eq:2d-signal-profile}
\begin{cases}
-\Delta v_k=am(1-\delta_k)(\phi_{\varepsilon_k}-\phi_\rho)
   &\text{in }\Omega,\\
\partial_\nu v_k=0 &\text{on }\partial\Omega.
\end{cases}
\end{equation}
Both rescaled profiles have integral one on $\Omega$, so the right-hand side
has zero integral and the Neumann problem is solvable.
Writing $P(s)=\int_{B_s}\phi$ for $s\ge0$, we integrate
\eqref{eq:2d-signal-profile} radially to obtain
\begin{equation}
\label{eq:2d-family-derivative}
2\pi r v_k'(r)=-am(1-\delta_k)
  \bigl[P(r/\varepsilon_k)-P(r/\rho)\bigr]
\quad\text{for }r\in(0,R],\ k\ge1.
\end{equation}
Since $\varepsilon_k\le\rho$ and $P$ is nondecreasing, $v_k'\le0$.
Thus $v_k\ge1$. The density $u_k$ is positive, has mass $m$,
and is constant near $\partial\Omega$, so $\partial_\nu u_k=0$
there. Both functions are smooth and radial, which proves
\eqref{eq:2d-family-mass-norm}. Also,
\eqref{eq:2d-family-derivative} yields
\begin{equation*}
1\le v_k(r)\le1+\frac{am}{2\pi}\ln\frac{R}{r}
\quad\text{for }r\in(0,R],\ k\ge1.
\end{equation*}
In particular, $\|v_k\|_2$ is bounded uniformly in $k$.

We compute the energy using \eqref{eq:2d-exact-energy-representation}.
By convexity and superadditivity of $s\mapsto s\ln s$, the entropy
of the sum in \eqref{eq:2d-concentrating-data} differs from the sum
of the two entropies by at most $m\ln2$. Rescaling the concentrated
term and using $\delta_k\ln(1/\varepsilon_k)=O(1)$ gives
\begin{equation}
\label{eq:2d-family-entropy}
\int_\Omega u_k\ln u_k=2m\ln\frac1{\varepsilon_k}+O(1)
\quad\text{as }k\to\infty.
\end{equation}
Let $M_k,E_k$ be the cumulative quantities in
\eqref{eq:2d-mass-definitions} for $(u_k,v_k)$. By
\eqref{eq:2d-radial-signal} and \eqref{eq:2d-family-derivative},
\begin{equation*}
E_k(r)=(1-a)M_k(r)+am(1-\delta_k)P(r/\rho)
                 +am\delta_k\frac{r^2}{R^2}
\quad\text{for }r\in[0,R],\ k\ge1.
\end{equation*}
Because $\rho$ is fixed and $P(s)\le C\min\{s^2,1\}$ for $s\ge0$,
the last two terms are bounded by $Cr^2$ uniformly in $k$.
Also, $M_k(r)\le C(r/\varepsilon_k)^2$ for $r\in[0,\varepsilon_k]$,
and $|M_k(r)-m|\le m\delta_k$ for $r\in[\varepsilon_k,R]$.
Therefore,
\begin{align}
\int_0^R\frac{E_k(r)^2-M_k(r)^2}{r}\dd r
&=-a(2-a)\int_0^R\frac{M_k(r)^2}{r}\dd r+O(1)\notag\\
&=-\frac{m^2b}{1+b}\ln\frac1{\varepsilon_k}+O(1)
\quad\text{as }k\to\infty.
\label{eq:2d-family-mass-energy}
\end{align}
Indeed, the contribution of $M_k^2/r$ from $(0,\varepsilon_k)$
is bounded, and the error from replacing $M_k^2$ by $m^2$ on
$[\varepsilon_k,R]$ is $O(\delta_k\ln(R/\varepsilon_k))=O(1)$.
Combining \eqref{eq:2d-family-entropy} and
\eqref{eq:2d-family-mass-energy} in \eqref{eq:2d-exact-energy-representation}
yields
\begin{align*}
\mathcal F_{\mathrm{KS}}(u_k,v_k)
&=\left[2m(1+b)-\frac{m^2b}{4\pi}\right]\ln k+O(1)\\
&=2m\left(1-\frac{m}{8\pi}\right)\ln k+O(1)
\quad\text{as }k\to\infty,
\end{align*}
where we used the uniform $L^2$ bound on $v_k$ and $v_k(R)=1$.
This proves \eqref{eq:2d-family-energy}.

For the dissipation, let $f_k,g_k$ be given by
\eqref{eq:2d-dissipation-variables} for $(u_k,v_k)$. By
\eqref{eq:2d-concentrating-data} and \eqref{eq:2d-signal-profile},
\begin{equation*}
f_k=-(1-a)m(1-\delta_k)\phi_{\varepsilon_k}
    +v_k-am(1-\delta_k)\phi_\rho-\frac{m\delta_k}{|\Omega|}
\quad\text{in }\Omega.
\end{equation*}
The last three terms are uniformly bounded in $L^2(\Omega)$, so
\begin{equation*}
\|f_k\|_2=(1-a)m(1-\delta_k)\varepsilon_k^{-1}\|\phi\|_{L^2(B_1)}+O(1)
\quad\text{as }k\to\infty.
\end{equation*}
The bound on $P$ and \eqref{eq:2d-family-derivative} give
$\|\nabla v_k\|_\infty\le C\varepsilon_k^{-1}$.
Since $u_k$ is a positive constant plus the square of
$\sqrt{m(1-\delta_k)}\,\varepsilon_k^{-1}\psi(x/\varepsilon_k)$,
\begin{equation*}
\int_\Omega\frac{|\nabla u_k|^2}{u_k}
\le \frac{4m(1-\delta_k)}{\varepsilon_k^2}
    \|\nabla\psi\|_{L^2(B_1)}^2
\quad\text{for }k\ge1.
\end{equation*}
Consequently,
$\|g_k\|_2^2\le 2\int_\Omega|\nabla u_k|^2/u_k
   +2m\|\nabla v_k\|_\infty^2\le C\varepsilon_k^{-2}$.
By \eqref{eq:2d-dissipation-decomposition} and $1-a>0$, we obtain
\begin{equation*}
c\varepsilon_k^{-2}
\le\mathcal D_{\mathrm{KS}}(u_k,v_k)
\le C\varepsilon_k^{-2}
\quad\text{for sufficiently large }k.
\end{equation*}
Taking logarithms yields \eqref{eq:2d-family-log-dissipation} because
$2(1+b)=4(m-4\pi)/(m-8\pi)$.
The ratio of $-\mathcal F_{\mathrm{KS}}(u_k,v_k)$ to
$\ln(1+\mathcal D_{\mathrm{KS}}(u_k,v_k))$ thus has a strictly
positive limit. This rules out replacing $\ln(1+s)$ by
$h(s)=o(\ln(1+s))$ in \eqref{eq:2d-sobolev-energy-dissipation}.
\end{proof}

\section*{Declarations}
\begin{itemize}
  \item Conflict of interest: The authors declare they have no financial interests.
  \item Data availability: Not applicable.
\end{itemize}

\subsection*{Declaration of generative AI use}
The authors used OpenAI Codex to assist with language editing and
\LaTeX{} preparation, to explore proof strategies, and to refine the
exposition of mathematical arguments. The authors take full
responsibility for the final content of the manuscript.

\raggedbottom
\providecommand{\bysame}{\leavevmode\hbox to3em{\hrulefill}\thinspace}
\providecommand{\MR}{\relax\ifhmode\unskip\space\fi MR }
\providecommand{\MRhref}[2]{%
  \href{http://www.ams.org/mathscinet-getitem?mr=#1}{#2}
}
\providecommand{\href}[2]{#2}

\end{document}